\documentclass{article}

\usepackage{arxiv}

\usepackage[utf8]{inputenc} 
\usepackage[T1]{fontenc}    
\usepackage{hyperref}       
\usepackage{url}            
\usepackage{booktabs}       
\usepackage{amsfonts}       
\usepackage{nicefrac}       
\usepackage{microtype}      
\usepackage{graphicx}
\graphicspath{ {./images/} }
\usepackage[normalem]{ulem}
\usepackage{xcolor}%
\usepackage[authoryear]{natbib}
\usepackage{pdflscape}

\usepackage[ruled,vlined]{algorithm2e}

\usepackage{libertine}
\usepackage{libertinust1math}
\usepackage{dsfont}  
\usepackage{nicefrac}
\usepackage{todonotes}

\usepackage{graphicx}
\graphicspath{{media/}}     
\usepackage{libertine}
\usepackage{libertinust1math}
\usepackage{dsfont}  
\usepackage{nicefrac}
\usepackage{todonotes}

\usepackage{float} 
\usepackage{lscape}

\title{On Solving Large Scale Dynamic Weapon Target Assignment Problems}

\author{
 Radosław Pytlak \\
  Faculty of Mathematics and Information Science\\
  Warsaw University of Technology\\
  Koszykowa 75, Warsaw, 00-662\\ Poland \\ \texttt{radoslaw.pytlak@pw.edu.pl}\\
   \And
 Robert Przybycień \\
  Faculty of Mathematics and Information Science\\
  Warsaw University of Technology\\ Koszykowa 75, Warsaw, 00-662\\ Poland\\ \texttt{robert.przybycien.dokt@pw.edu.pl}\\
  \And
 Kamil Kisiel \\
  Faculty of Mathematics and Information Science\\
  Warsaw University of Technology\\ Koszykowa 75, Warsaw, 00-662\\ Poland\\
\texttt{kamil.kisiel3.dokt@pw.edu.pl} \\
  \And
 Wojciech Stecz \\
  Institute of Computer Sciences\\ Cardinal Stefan Wyszyński University\\ Wóycickiego 21, Warsaw, 01-938\\ Poland\\ \texttt{w.stecz@uksw.edu.pl}\\  
  }

\begin{document}
\maketitle
\begin{abstract}
The problem of Weapon Target Assignments (WTA) has long attracted researchers' interest due to its practical significance and the fact that it is an NP--complete, nonlinear convex problem. Significant progress has been made in recent years regarding algorithms and solvers for this problem. Recent approaches to solving the WTA problem take into account the convexity of the objective function, employ problem reformulation to facilitate the efficient use of column generation, and propose a hybrid approach---in this method, solving a relaxed convex problem in the initial stage provides a starting set of columns for the column generation process---enabling the efficient solution of static WTA problems involving more than 10 000 targets and weapons. Practical considerations, however, necessitate defining the WTA problem under the assumption that new targets appear randomly over time. This leads to the formulation of the WTA problem as a multistage problem. In this paper, we present a formulation of the dynamic WTA problem and demonstrate, based on numerous numerical experiments, that it can be effectively solved using a strategy rooted in a scalable Branch-Price-and-Cut algorithm, employing a receding horizon control technique.
\end{abstract}

\section{Introduction}
The problem of optimally assigning destructive means to air targets that threaten protected facilities is called the Weapon Target Assignment (WTA) problem. This problem is one of four looping activities in the process of combating air threats in a cycle called "observe-orient-decide-act" (OODA) \cite{Zhao2022}. OODA is an operational framework developed by United States Air Force Colonel John Boyd in the 1970s. It consists of four sequential tasks: Threat Assessment, Target Selection, Weapon Target Assignment, and Efficiency Evaluation.
This article focuses on methods for solving the third WTA task. The preceding tasks allow assessment of the threat posed by approaching strike assets and preliminary target selection. 

The Threat Assessment process makes it possible to determine the probability of specific aerial objects striking protected targets. This is often accomplished by combining models that classify the type of airborne asset (UAV, fighter, missile, etc.), providing insight into parameters such as range and effectiveness, and models that assess the target's intent based on observed behavioral characteristics---such as distance from the defended zone, history, and trajectory. These methods rely on algorithmic inference from the provided potential threat parameters, representing knowledge as ontologies or Bayesian networks (probabilistic graph models---see Johansson et al. \cite{4632368}). In practice, to assess the feasibility of attacking selected targets, heuristic algorithms determine the flight path of aerial assets, accounting for the aircraft's aerodynamic capabilities, the terrain grid, and the detection capabilities of defensive systems. One frequently used efficient algorithm is the Rapidly-exploring Random Tree \cite{burzynski2024trajectory}.

The Target Selection task selects targets to be destroyed first based on the results of the Threat Assessment task. This task stems from the assumption that, in some cases, the number of incoming targets exceeds the capabilities of air defenses, so the most significant targets should be selected. These are typically the targets that can pose the greatest threat.

The WTA task has two basic variants: Static WTA and Dynamic WTA \cite{LI2024109212, KLINE2019226}. The static task assigns destructive agents to threats without considering when threats will be within range of the destructive agents. In this case, we assume that any threat can be destroyed. Dynamic WTA is an extended task that assigns destructive agents to targets in subsequent time intervals. These intervals depend on numerous factors, may include the moments of
launching destructive agents at targets, the moments of hits or misses, etc. The detailed process of generating these intervals is beyond the scope of this article. However, we assume we can uniquely determine these intervals, which is usually relatively simple. 

The final task of the OODA loop (Efficiency Evaluation) assesses achieved targets. Based on the verified results of destroying incoming targets, the system moves on to the next threat assessment in the OODA loop. The presented principle of combating threats involves executing individual tasks in a loop. This model is often referred to as shoot-look-shoot. 

The WTA problem itself is not a new optimization problem. WTA belongs to a family of nonlinear assignment problems, including facility location, utility allocation, and radiation treatment problems. This means solution approaches developed and implemented in one application domain can be applied to other problems in this family.

The following are contributions of our paper.
\begin{enumerate}
    \item Proposing a method for solving the static WTA problem using the objective function linearization scheme proposed by Andersen et al. \cite{andersen} and the Branch--and--Price approach.
    
    \item Application of a hybrid approach combining the Andersen \cite{andersen} and Lu-Chen algorithms with the initialization proposed by Bertsimas to solve the static WTA problem.

    \item Numerical comparison of the Andersen et al. \cite{andersen}. and Le at al. \cite{chen} methods—equipped with the Bertsimas et al. \cite{bertsimas} initialization procedure---with the Bertsimas et al. \cite{bertsimas} approach.

    \item Proposing a method for solving the dynamic WTA problem based on the Model Predictive Control approach.
  
\end{enumerate}

The paper is organized as follows. Section \ref{sec:RelatedWorks} presents a review of the literature required to evaluate the proposed approach to solving the dynamic WTA problem.. Section \ref{sec::Methods} describes WTA problems, divided into static problem presented in Section \ref{sec::StaticWTA} and dynamic problem presented in Section \ref{sec::DynamicWTA}. The results obtained for both problem types are presented in Section \ref{sec:Results}. Section \ref{sec:Discussion} analyzes the practical application of the described algorithms. Further directions of algorithm development are presented in Section \ref{sec:Conclusions}.

\section{Related Works} \label{sec:RelatedWorks}
In recent years, a number of optimization methods for static and dynamic WTA problems have been proposed. Among the proposed solution methods, heuristic algorithms are the most frequently employed, although optimization methods for convex integer programming problems utilizing column generation algorithms also play a significant role. Readers interested in this topic may consult two review articles \cite{LI2024109212, KLINE2019226}, which present new algorithms for solving static and dynamic versions of the WTA problem. In this paper, we compare our proposed methods with selected approaches described in the aforementioned publications. The comparisons we present focus primarily on exact methods based on efficient integer programming techniques designed for static WTA problems: the method of Andersen et al. \cite{andersen}, the algorithm of Lu et al. \cite{chen}, and the hybrid approach proposed by Bertsimas et al. \cite{bertsimas}.

Leading methods employing exact integer optimization algorithms will be discussed in great detail in subsequent sections of this paper; in this part, we focus on heuristic and metaheuristic methods. 
Examples of using heuristics to solve WTA appear in the papers by Hendrickson et al. \cite{ hendrickson2023decentralized} and Ma et al. \cite{ MA2021107717}. 
Hendrickson et al. \cite{ hendrickson2023decentralized} introduce a method for distributed assignment planning executed by the weapons, where each weapon solves an optimization problem over distinct subsets of the decision variables. 
Ma et al. \cite{ MA2021107717} present a WTA problem that accounts for uncertainty in assigning weapons to threats. To address the uncertainty of weapon adequacy, the authors describe the WTAU (WTA with uncertainty) problem as a nonlinear integer programming model. To improve computational efficiency, they present a hybrid heuristic search algorithm. Dagistanli et al. \cite{DAGISTANLI2026115932} proposed an approach to solving this problem based on selected heuristic methods. 
There are a number of heuristic methods for solving static WTA problems that are based on metaheuristics such as: MOEA/D (Multi-Objective Evolutionary Algorithm based on Decomposition), and NSGA-II (Non-dominated Sorting Genetic Algorithm II), Artificial Fish Swarm Algorithm (AFSA), or PSO (Particle Swarm Optimization). 
Xu et al. \cite{XU2020112844} formulated a two--objective, multi--stage WTA task allocation model that analyzes interactions between sensor and weapon platforms. A key innovation is using the Soyster model's robust parameters to address uncertainty in real-time task allocation. To ensure flexibility, the authors employ a decomposition--based MOEA/D.
Gong et al. \cite{10.1007/978-981-99-9021-4_58} solve a multi--objective WTA problem with NSGA-II. This method demonstrates NSGA-II's feasibility for solving this dynamic WTA problem across three datasets. Xu et al. \cite{ XU20242694} propose a constrained optimization model, Firepower Conflict Free WTA (FCFWTA), and a Survival Evolution (SE) strategy for the AFSA, this method solves the complex constrained WTA problem. 

In the case of dynamic WTA problem the theoretically justified approach based on Bellman's equations is presented in Section \ref{sec::DynamicWTA}, in this section we mention two basic heuristic approaches to dynamic WTA. In the first approach the classical dynamic WTA is considered in which we have horizon with $N$ stages and on which stage weapons should be allocated in such a way that the targets destruction is evaluated over whole horizon---\cite{Xin2010}, \cite{Xin2011}, \cite{Chang1987}, \cite{Leboucher2013}, \cite{Karasakal}, \cite{Arslan}---in \cite{Khosla2001}, \cite{Wu2008}, \cite{Chen2009}, \cite{CHANG2018845}, \cite{XU2020112844}, \cite{Kong}, methods based on genetic algorithms have also been proposed. A typical example of a multi-stage dynamic WTA problem is the work of Chan et al. \cite{CHANG2018845} and Pardalos et al. \cite{XU2020112844}. Chang et al. described a method for constructing a multi-stage plan and the shoot--look--shoot principle. A heuristic method called the improved artificial bee colony algorithm with heuristic factor initialization was used to solve the dynamic WTA problem. Pardalos et al. presented a heuristic approach to solving a multi-stage problem based on the MOEA/D evolutionary algorithm. Kong et al. \cite{Kong} solve the same problem using a PSO heuristic. 

The second approach involves solving the dynamic WTA problem in two steps. In the first step, available resources are allocated to incoming threats. In the second step, static WTA problems are solved for all stages---Zhao et al. \cite{Zhao2022} presented this approach. However, planning resource allocation for individual stages requires precise knowledge of the specifics of incoming air threats. It also depends heavily on the type of defense means available to the defending side. However, this approach is close to real-world operations. In \cite{8668760} an et al. present a different approach, the authors model the dynamic WTA problem using game theory. However, this approach requires an extensive mathematical apparatus, and solving problems formulated within this framework necessitates making highly simplifying assumptions which are far from the principles of planning actions in real-world conditions.

With the increased computational capabilities of air defense systems and the use of efficient sensors that can track multiple targets simultaneously, researchers are developing Sensor Weapon Target Assignment (SWTA) models. Wang et al. \cite{10621019} use the Large Neighborhood Search heuristic to solve the SWTA problem. Zhang et al. \cite{9288732} present the EA-DS/HWTA evolutionary algorithm to solve the SWTA problem. Xin et al. \cite{8249748} proposed the MRBCH heuristic to solve the formulated sensor-weapon-target assignment (S-WTA) problem.

\section{Methods} \label{sec::Methods}

\subsection{Static WTA problem} \label{sec::StaticWTA}

The basic WTA problem aims to allocate weapons to targets to minimize the targets 'expected survival value. The problem is formulated as follows. The objective function minimizes the targets' survival value:

\begin{equation}
\min_{x_{ij}} \sum\limits_{j = 1}^N V_{j} \cdot \prod\limits_{i = 1}^M (1-p_{ij})^{x_{ij}}
\label{eq:SWTAfcelu}
\end{equation} 
subject to the constraints:
\begin{equation}
\sum\limits_{j = 1}^N x_{ij} =1, \forall i = 1,\dots,M
\label{eq:SWTAfogr_a}
\end{equation} 
\begin{equation}
x_{ij} \in \{0,1\}, \forall j = 1,\dots,N, i =1,\ldots,M
\label{eq:SWTAfogr}
\end{equation} 

The variable $x_{ij}$ is set to 1 when launcher $i$ is assigned to strike target $j$, otherwise $x_{ij}$ is set to 0. 
The parameter $V_j$ specifies the target value, which allows them to be prioritized to determine the order in which they are destroyed. $M$ is the number of launchers and $N$ is the number of targets.

The WTA problem (\ref{eq:SWTAfcelu})--(\ref{eq:SWTAfogr}) is a nonlinear integer programming problem in which the solution scales nonlinearly with the number of weapons assigned to targets. The main challenge in solving the WTA problem is handling the nonlinear objective function.
In practice, solving a $200\times 200$ problem (number of weapons vs. number of targets) is not possible in less than several tens of minutes. 

To apply exact integer optimization methods to solve the problem, it is necessary to reformulate it. Andersen et al. \cite{andersen} proposed a formulation as a linear mixed integer optimization problem in which the objective function is linearized and linearization parameters are continuous decision variables.
Lu et al. \cite{chen} and Bertsimas et al. \cite{bertsimas} propose the reformulation based on Dantzig--Wolfe decomposition (\cite{Barnhart}). Although this approach involves a significantly larger number of decision variables than the original formulation, the use of efficient enumeration or column generation procedures mitigates this issue. Both reformulations are presented in considerable detail later in the paper..

\subsection{Linearization of the basic WTA problem}
\label{sec::WTAlinearization}
Linearizing the objective function, as proposed in \cite{andersen}, allows us to prepare a linear formulation of the WTA problem while providing a good approximation of the objective function (\ref{eq:SWTAfcelu}). We rely on the approach presented by Camm et al. \cite{Camm} and Andersen et al. \cite{andersen}. 

Consider a specific target $j$, then the probability of its survival is evaluated as:

\begin{equation}
z_{j} = \prod\limits_{i = 1}^M (1-p_{ij})^{x_{ij}} \label{eq::SWTAprof}
\end{equation}

In order to avoid using the product (\ref{eq::SWTAprof}) in the objective function we take the logarithm of variable $z_i$ to obtain:
\begin{equation}
\ln(z_j) = \sum\limits_{i = 1}^M x_{ij} \cdot \ln(1-p_{ij})
\label{zj}
\end{equation} 

The problem (\ref{eq:SWTAfcelu})-(\ref{eq:SWTAfogr}) is now transformed into the problem 

\begin{equation}
\min_{x_{ij}} \sum\limits_{j = 1}^N V_{j} \cdot e^{\ln(z_j)}
\label{eq:SWTAlinfcelu}
\end{equation} 
subject to the constraints:
\begin{equation}
\ln(z_j) = \sum\limits_{i = 1}^M x_{ij} \cdot \ln(1-p_{ij}),\  \forall j = 1,\dots,N
\label{zja}
\end{equation} 
\begin{equation}
\sum\limits_{j = 1}^N x_{ij} = 1,\ \forall i = 1,\ldots,M \label{SWTAogrW}
\end{equation} 
\begin{equation}
x_{ij} \in \{0,1\},\ \forall j = 1..N,\ i = 1,\dots,M
\label{eq:SWTAlinfogr}
\end{equation} 

\begin{figure}[ht]
\centering
\includegraphics[width=0.7\linewidth]{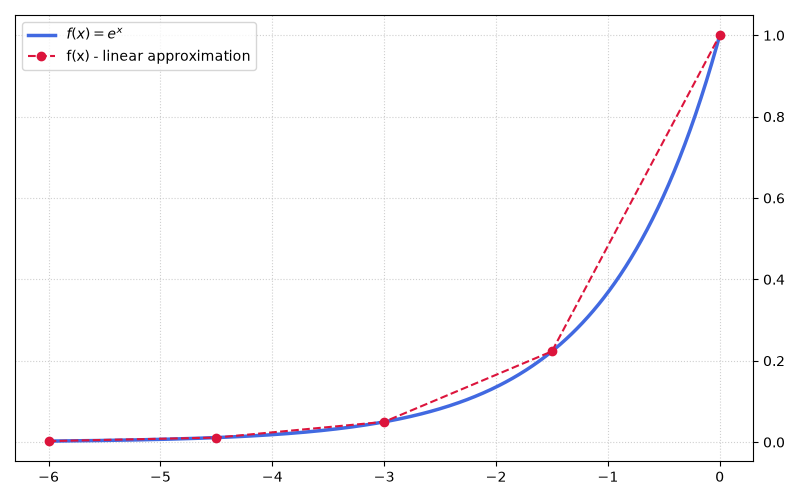} 
\caption{Example of linearization of a function $e^x$ using piecewise linear functions and five nodes $b_\gamma$, $\gamma \in \Gamma$ ($|\Gamma |=5$).} 
\label{fig::explinear}
\end{figure}

In order to linearize the above problem, it is assumed that the nonlinear element $e^{\ln(z_j)}$ will be approximated by a piecewise linear function as shown in Fig.\ref{fig::explinear}.

So, we assume that $f:{\cal R}\rightarrow {\cal R}$ is a convex piecewise linear function that approximates $e^x$ from above, and instead of solving the problem (\ref{eq:SWTAlinfcelu})-(\ref{eq:SWTAlinfogr}) we tackle the problem
\begin{equation}
\min_{x_{ij}} \sum\limits_{j = 1}^N V_{j} \cdot f(\ln(z_j))
\label{eq:SWTAlinfcelu_e}
\end{equation} 
subject to the constraints (\ref{zja})--(\ref{eq:SWTAlinfogr}).

It is shown in \cite{andersen} that the function $f$ can be found by solving the following problem
\begin{equation}
f(x) = \min_{\kappa_\gamma} \sum\limits_{\gamma \in \Gamma} \kappa_{\gamma} \cdot e^{b_\gamma}
\end{equation} 
subject to the constraints
\begin{equation}
\sum\limits_{\gamma \in \Gamma} \kappa_{\gamma} \cdot b_\gamma = x,
\end{equation} 
\begin{equation}
\sum\limits_{\gamma \in \Gamma} \kappa_{\gamma} = 1,
\end{equation} 
\begin{equation}
\kappa_{\gamma} \geq 0,\ \forall \gamma \in \Gamma
\end{equation} 

Once the breakpoints $b_{\gamma}$ are defined for every target in the WTA problem the linear mixed integer optimization problem can be defined in which $x_{ij}$ and $\lambda_{\gamma}$ are decision variables. The details of the procedure for choosing the breakpoints $b_{\gamma}$ and  the final form of the optimization problem are stated in \cite{andersen}. 
The number of the problem variables increases with the increased number of breakpoints $b_\gamma$,  notably, for the WTA problem with an exponential function, the number of nodes is small and limited to a few.

\subsection{Column generation algorithm for WTA} \label{sec::WTAGGdescription}

Column Generation is an iterative optimization method that solves large linear programs by splitting the problem into two components: a Restricted Master Problem (RMP) and a Pricing Subproblem (SUB). It starts with a small subset of columns (an initial set of patterns from which it chooses a feasible solution), uses dual variables to find new columns with negative reduced costs (associated with the optimal value of the SUB), and loops until no improving patterns remain. Desrosiers et al. provide a good introduction to the topic \cite{Desrosiers2005, Desrosiers2026}, Barnhardt et al. provide an advanced theoretical analysis \cite{Barnhart}.

\begin{figure}[ht]
\includegraphics[width=1.0\linewidth]{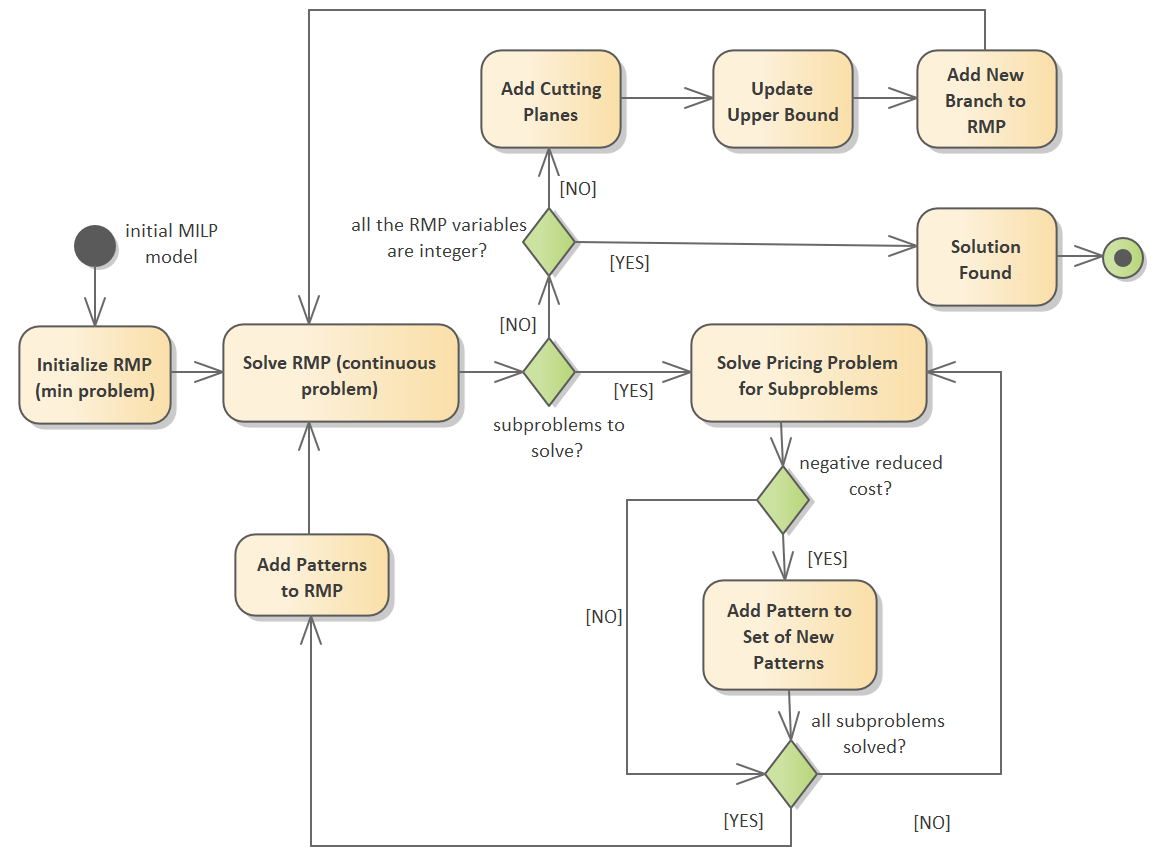} 
\caption{A scheme of the modified two--loop column generation algorithm. The loop containing the Solve RMP problem and the Solve Pricing Problem for Subproblems represent the basic problems of column generation and adding new columns/patterns to the RMP problem. The outer loop containing the problems Add Cutting Planes, Update Upper Bound, and Add New Branch to the RMP is executed when no additional patterns can be found for the given RMP form. However, the solution to the RMP problem contains variables whose values are not integers.} 
\label{fig::CGscheme}
\end{figure}

The method's popularity stems from the range of problems it can effectively solve. These include Vehicle Routing Planning (VRP), Task Assignment, and Production Planning (PP). A characteristic feature of these problems is that they can be divided into subproblems for which most constraints are defined in the definition of the optimization probleml. The remaining constraints, in significantly smaller numbers, connect the individual component problems in the RMP. In the case of the WTA, the subproblem is determining a list of weapons for an incoming target. The upper RMP problem selects those patterns from among potential threat destructive patterns that satisfy the connecting constraints. In the case of the WTA, this requires that the weapon used to destroy the target be used at most $w_i$ times (usually once). Fig.\ref {fig::CGscheme} shows a sketch of the algorithm.

The algorithm presented in Fig.\ref{fig::CGscheme} consists of several steps and two loops. The actions in the inner loop aim to find a feasible solution (pattern) for the given RMP that cannot be improved further; i.e., the reduced costs are negative (or positive, depending on the direction of optimization of the objective function).

First, the algorithm must determine a set of initial patterns from which it will select a feasible solution. For certain problems, finding a good set of initial patterns can be time-consuming. On the other hand, initial patterns that lead to a solution far from optimal can significantly slow the process of finding the optimal solution. This will be discussed in the article when presenting initialization heuristics for the WTA algorithm.

To efficiently find further feasible solutions, the RMP must be a linear problem. After the RMP determines the optimal relaxed solution (decision variables are continuous variables)  for a given set of patterns, it verifies whether to add new patterns by solving a valuation subproblem and checking the reduced costs. For this purpose, the values of dual variables of the RMP are used as parameters in the subproblems. For each subproblem, the possibility of adding new patterns using valuation is verified. In the case of WTA, the subproblem involves allocating weapons to each incoming target. If an allocation with negative reduced costs can be found, such a pattern is added to the set of patterns.

If patterns with reduced costs cannot be found, it is necessary to check whether all decision variables that are supposed to be integers are actually integers. If this is not the case, fractional  variables should be selected and their values determined using the Branch--and--Bound technique. This requires adding new constraints to the RMP problem on variables that must be checked for integer values (the outer loop of the algorithm). After adding new constraints, we solve the RMP and its subproblems again. The number of such iterations decreases with each new variable value constraint added to the RMP. Nevertheless, the branching method must be carefully defined to obtain fast results. We discuss this issue later in the paper.

A significant speedup in determining the solution to the problem can be achieved by introducing additional hyperplane cuts for a specific WTA problem, then we refer to the method as Branch--Price--and--Cut approach (\cite{bertsimas}). 

\subsection{Linear version of the WTA Branch-and-Price subproblem} \label{sec::WTAGGlinearization}
Using the WTA problem (\ref{eq:SWTAfcelu})-(\ref{eq:SWTAfogr}) to construct exact algorithms has not yielded efficient methods for the problem. This is because such algorithms rely on a relaxed version of the problem (\ref{eq:SWTAfcelu})-(\ref{eq:SWTAfogr}) (with variables $x_{ij} \in \{0,1\}$) that provides an imprecise lower bound on the optimal objective function value. Consequently, applying the Branch--and--Bound method results in highly complex branching trees when seeking high-quality feasible solutions. 

In order to overcome the problem with the low--quality of the optimal value estimation based on the relaxed problem Lu et al \cite{chen} and Bertsimas et al. \cite{bertsimas} use the reformulation of the WTA problem which has many more decision variables but at the same its relaxed problem guarantees much better estimation of the optimal value of the problem. 

To present the new formulation of the WTA problem, let us consider a given target $j$. There are $2^M - 1$ possible attack variants against this target---weapon $i$ may or may not be assigned to target $j$ (assuming no target is left unaddressed). We are thus dealing with all possible subsets of weapons used to engage the threats. Each such subset used to destroy threat $j$ can be viewed as an attack variant (pattern) against that threat. Let us assume there are $L$ such patterns. The new formulation of the WTA problem is then as follows:

\begin{equation}
\min_{y_{j}^l} \sum\limits_{j = 1}^N \sum\limits_{l = 1}^L q_{j}^l \cdot y_{j}^l
\label{eq::WTABBfcelu}
\end{equation} 
subject to the constraints:
\begin{equation}
\sum\limits_{j = 1}^N \sum\limits_{l = 1}^L b_{ij}^l \cdot y_{j}^l = 1,\ \forall i = 1,\ldots,M
\label{eq::WTABBconstraint1}
\end{equation} 
\begin{equation}
\sum\limits_{l = 1}^L y_{j}^l = 1,\ \forall j = 1,\dots,N
\label{eq::WTABBconstraint2}
\end{equation} 
\begin{equation}
y_{j}^l \in \{0,1\},\ \forall j = 1,\ldots,N, l = 1,\ldots,L
\label{eq::WTABBconstraint3}
\end{equation} 

\begin{figure}[ht]
\centering
\includegraphics[width=0.7\linewidth]{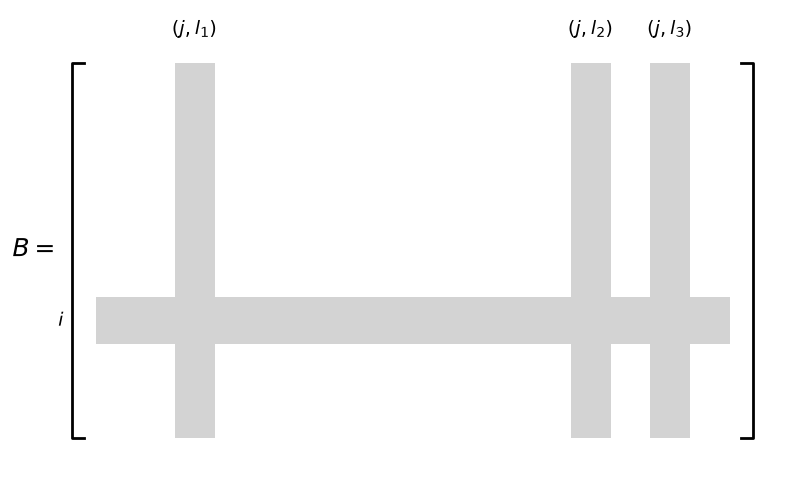} 
\caption{RMP matrix. Matrix B is a matrix containing in each column a separate pattern of destruction of the $j$th target by the $i$th weapon. Solving the column generation subproblem means adding a new pattern $x_j$ as a column $B_j^l$ of matrix $B$. Three target destruction patterns for $j$th threat are highlighted.} 
\label{fig::RMPmatrixEx}
\end{figure}

(\ref{eq::WTABBfcelu}) states that we minimize the value of survived targets. The variable $y_{j}^l$ equals 1 when pattern $l$ for using a subset of weapons against target $j$ is selected. The parameter $q_{j}^l$ represents the target value after firing under scenario representing by pattern $l$. The constraint (\ref{eq::WTABBconstraint1}) says that each weapon can be used only once. Parameter $b_{ij}^l$ equals 1 if weapon $i$ is used to destroy target $j$ in pattern $l$. The constraint (\ref{eq::WTABBconstraint2}) forces the selection of one scenario for each threat.

In the RMP model of the problem (\ref{eq::WTABBfcelu})--(\ref{eq::WTABBconstraint3}) $q_{j}^l$ is a parameter and is evaluated as follows (see goal function (\ref{eq:SWTAfcelu})):
\begin{equation}
q_{j}^l = V_j \cdot \prod\limits_{i = 1}^M [1 - b_{ij}^l \cdot p_{ij}] = V_j \cdot \prod\limits_{i = 1}^M [1 - p_{ij}]^{b_{ij}^l}
\label{qjsdef4}
\end{equation}

In a column generation procedure we have to formulate the subproblem which solution defines a new pattern. To this end we refer to the problem of the relaxed RMP (\ref{eq::WTABBfcelu})--(\ref{eq::WTABBconstraint3}). Suppose that $\mu \in {\cal R}^N$ are the Lagrange multipliers of the problem associated with the constraints (\ref{eq::WTABBconstraint2}) and $\lambda \in {\cal R}^M$ are the Lagrange multipliers associated with the constraints (\ref{eq::WTABBconstraint1}). The reduced cost of the pattern $l$ (see for example \cite{bertsimas}) evaluated for the target $j$ is expressed as:
\begin{equation}
    q_j^l - \mu_j - \sum_{i=1}^M b^l_{ij} \lambda_i \label{redCost}
\end{equation}
If the reduced cost, for a given pattern, is negative then incorporating this pattern to the current set of patterns can reduce value of the objective function of the relaxed RMP. Therefore, for each target $j=1,\ldots,N$ we solve the subproblem which we call the Pricing Problem
\begin{equation}
    \min_{x\in \{0,1\}^M} \left [ V_j \prod_{i=1}^M (1-p_{ij})^{x_{ij}} -\lambda^Tx - \mu_j\right ]. \label{PP}
\end{equation}

Pricing is a decision step in the Simplex algorithm. It involves calculating the reduced costs of non-basic variables to select the one that will most rapidly improve the objective function value. In the context of the column generation method, it allows for verifying whether a determined column (a new pattern) reduces the value of the objective function.

By taking into account (\ref{eq::SWTAprof})--(\ref{zj}) the above problem can be restated as
\begin{equation}
    \min_{x\in \{0,1\}^M} \left [ V_j e^{\ln (z_j)}  -\lambda^Tx - \mu_j\right ] \label{PPa}
\end{equation}
subject to the constraints
\begin{equation}
\ln(z_j) = \sum\limits_{i = 1}^M x_{ij} \cdot \ln(1-p_{ij}).
\label{PPb}
\end{equation}

Furthermore, if we introduce vector $s_j = [s_{1,j},s_{2,j},\ldots,s_{M,j}]$, where $s_{i,j} = \ln (1-p_{ij})$ then the problem (\ref{PPa}) can also be stated in the following way
\begin{equation}
    \min_{x\in \{0,1\}^M} \left [ V_j e^{s_j^Tx} -\lambda^Tx - \mu_j\right ]. \label{PPc}
\end{equation}

We have at least two approaches to solving the Pricing Problem. The first one is based on the linearization of the problem (\ref{PPa})--(\ref{PPb}) using the procedure introduced by Andersen et al. \cite{andersen}. The second approach refer to the relaxation of the problem (\ref{PPc}) which will the convex optimization problem with constraints $0\leq x_{ij} \leq 1$, $i=1,\ldots,N$. The second approach is discussed in Subsection \ref{sec::WTAGGBertsimas}.

In the first approach we use the same linearization method as described in the Subsection \ref{sec::WTAlinearization}. It means that the function $q_j(x_{ij}$ is approximated by using breakpoints and numbers $\kappa_{j\gamma}$.
Once the breakpoints for the $j$th subproblem are determined the following mixed integer programming problem is solved 
\begin{equation}
    \min_{x\in \{0,1\}^M,\\
    \{\kappa_{j\gamma}\in {\cal R}\}_{\gamma \in \Gamma}} \left [V_j \cdot \sum\limits_{\gamma \in \Gamma} \kappa_{j\gamma} \cdot e^{b_{j\gamma}} -  \sum\limits_{i = 1}^M \lambda_i \cdot x_{ij} - \mu_j \right ] 
\end{equation}

subject to:
\begin{equation}
\sum\limits_{\gamma \in \Gamma} \kappa_{j\gamma} \cdot b_{j\gamma} = \sum\limits_{i = 1}^M x_{ij} \cdot \ln(1-p_{ij}) 
\end{equation}
\begin{equation}
\sum\limits_{\gamma \in \Gamma} \kappa_{j\gamma} = 1,\ \kappa_{j\gamma} \geq 0,\ \gamma \in \Gamma,
\end{equation}
where $|\Gamma|$ is equal to the number of used breakpoints in the approximation.

\subsection{An approach to solving the pricing subproblem by enumeration.} \label{sec::WTAGGChen}

The linear formulation (\ref{eq::WTABBfcelu})--(\ref{eq::WTABBconstraint3})  allows for another method of adding successive patterns $y_{js}$ without using linearization in the subproblems. According to \cite{chen}, standard pricing is difficult to implement due to the nonlinear dependence of $q_{j}^l$ on $x_{ij}$, $j=1,\ldots,M$. They propose to circumvent this problem by enumerating successive patterns for a given target $j$ as bit sequences of pattern $l$. These are then iterated until one with a negative pricing is found. The exponentially large number of bit sequences to be checked can be effectively reduced by using \textit{"weapon number bounding"}, which sets lower and upper limits of $1$ in the sequence, and \textit{"weapon domination"}, which, by understanding which weapons $i$ are always more effective for a given target, allows for bypassing some of the sequences being checked.

\subsection{WTA Branch-Price-and-Cut Linear Model} \label{sec::WTAGGBertsimas}

A fundamentally new approach to solving the static WTA problem was presented by Bertsimas et al. \cite{bertsimas}. This approach incorporates two elements that make it significantly more effective than those presented in the works of Andersen et al. \cite{andersen} and Lu et al. \cite{chen}. First, the new method represents a hybrid approach to solving a nonlinear integer optimization problem.In the first step of the method, a continuous nonlinear optimization problem is solved to determine initial patterns for the problem (\ref{eq::WTABBfcelu})--\ref{eq::WTABBconstraint3}). Subsequently, a static WTA problem based on the problem (\ref{eq::WTABBfcelu})--(\ref{eq::WTABBconstraint3}) is solved, with the pricing subproblem for target $j$ being solved as a continuous convex optimization problem subject to the constraints $0 \leq x_{ij} \leq 1$, $i=1, \ldots, M$, after which the non--integer values are converted to $0$ or $1$.

\subsubsection{Block coordinate descent algorithm for the Initialization}
Let $s_{ij}$ be a real-valued parameter associated with the survival probability of target $j$ as a result of employing weapon $i$: $s_{ij} = \ln \left(1 - p_{ij}  \right)$, then
the initialization problem can be stated as: 
\begin{equation}
\min_{X \in {\left\{ 1,0 \right\}}^{M\times N}} \sum\limits_{j = 1}^N V_j e^{s_j^T x_j}\label{InitPa}
\end{equation}
where $s_j = \left[s_{1j}, s_{2j},\ldots,s_{Mj}\right ]$. Since this is nonlinear integer programming problem Bertsimas et al. \cite{bertsimas} propose to solve the relaxed version of the problem (\ref{InitPa}): 
\begin{equation}
\min_{X \in {\cal R}^{M\times N}} \sum\limits_{j = 1}^N V_j e^{s_j^T x_j} \label{InitPb}
\end{equation}
\begin{equation}
\sum\limits_{j = 1}^N x_{ij} = \alpha \quad \forall i \in M, \label{InitPc}
\end{equation}
\begin{equation}
    0\leq x_{ij} \leq 1,\ i=1,2,\ldots,M,\ j=1,\ldots,N.
\end{equation}
where $x_j = \left [x_{1j},x_{2j},\ldots,x_{Mj}\right ]$. Notice that the relaxation is applied not only to binary variables $x_{ij}$ but also with respect to the constraints (\ref{SWTAogrW}).

In order to solve the problem (\ref{InitPb})--(\ref{InitPc}) Bertsimas et al. \cite{bertsimas} reformulate the problem to the following one
\begin{equation}
    \min_{X\in [0,1]^{M\times N}} \left [ \sum\limits_{j = 1}^N V_j e^{s_j^T x_j}   + \sum_{i=1}^M {\rm ind}_{\infty} \left (\sum_{j=1}^N x_{ij} =\alpha \right )\right ] \label{InitPd}
\end{equation}
where the indicator ${\rm ind}_{\infty}(\cdot)$ assumes value $+\infty$ when $\sum_{j=1}^N x_{ij} \neq\alpha $ and $0$ otherwise.

The indicator terms are separable with respect to the rows of the matrix $X$ thus the problem (\ref{InitPd}) can be stated as $N$ independent problems indexed by $i=1,2,\ldots,M$:
\begin{equation}
    \min_{x^i \in [0,1]^M} \sum_{j=1}^N U_j^i e^{s_{ij}x_{ij}} \label{InitPe}
\end{equation}
subject to
\begin{eqnarray}
    \sum_{j=1}^N x_{ij} = \alpha,\label{InitPf}
\end{eqnarray}
where $x^i = \left [ x_{i1},x_{i2},\ldots,x_{iM}\right ]$ and $U_j^i = V_j \sum_{k=1,\ldots,M,k\neq i} e^{x_{kj}x_{kj}}$.

Each of the problems (\ref{InitPe})–(\ref{InitPf}) is a constrained convex optimization problem satisfying Slater's condition; consequently, the Karush-Kuhn-Tucker necessary optimality conditions can be applied to solve them. Fix $i$, assume that $\lambda$ is a Lagrange multiplier corresponding to the constraint (\ref{InitPf}), $\gamma_j$ are the Lagrange multitpliers corresponding to the constraints $-x^i \leq 0$ and $\sigma_j$ corresponding to the constraints $x^i - 1 \leq 0$ then these optimality conditions are
\begin{equation}
    U_j^i s_{ij} e^{s_{ij}x_{ij}} + \lambda -\gamma_j + \sigma_j =0,\ \sum_{j=1}^N x_{ij} = \alpha,\ -x_{ij} \leq 0,\ x_{ij}-1 \leq 0,\ j=1,\ldots,N \label{InitPg}  
\end{equation}
\begin{equation}
    x_{ij} \gamma_{ij} =0,\ (x_{ij}-1)\sigma_j =0,\ \gamma_{ij} \geq 0,\ \sigma_j \geq 0,\ j=1,\ldots,N. \label{InitPh}
\end{equation}

In order to find $x_{ij}$, $j=1,2,\ldots,N$ satisfying the conditions (\ref{InitPg})--(\ref{InitPh}) we consider three cases:
\begin{enumerate}
\item If $\lambda \geq -u_js_{ij}$ then $\gamma_j - \sigma_j > u_j s_{ij}(e^{s_{ij} \hat x_{ij}} - 1) \geq 0$, and $\gamma_j > \sigma_j \geq 0$, so eventually $\hat x_{ij} = 0$
\item If $\lambda \leq -u_js_{ij} e^{s_{ij}}$ then $\gamma_j - \sigma_j < u_j s_{ij}(e^{s_{ij} \hat x_{ij}} - e^{s_{ij}}) \leq 1$, and $\sigma_j > \gamma_j \geq 0$, so eventually $\hat x_{ij} = 1$ 
\item If $-U_j^is_{ij} e^{s_{ij}} < \lambda < -U_j^is_{ij}$ then using previous arguments we can show that $U_j^i s_{ij}(e^{s_{ij} \hat x_{ij}} - e^{s_{ij}}) < \gamma_j - \sigma_j < U_j^i s_{ij}(e^{s_{ij} \hat x_{ij}} - 1)$. Assuming $\gamma_j > 0$ one obtains $x_{ij} = 0$ resulting in $\sigma_j > \gamma_j > 0$, which is a contradiction because $\gamma_j$ and $\sigma_j$ cannot both be non--zero. Similarly, if $\sigma_j > 0$ then $x_{ij} = 1$ and $\gamma_j > \sigma_j > 0$. Therefore, we should assume $\sigma_j = \gamma_j = 0$ resulting in $x_{ij} = \frac{1}{s_{ij}} \ln \left( -\frac{\lambda}{U_j^i s_{ij}} \right)$ 
\end{enumerate}
An analysis of the three cases leads to the conclusion that determining the decision variables and Lagrange multipliers that satisfy the KKT conditions reduces to finding a value of $\lambda$ for which the equation $\sum_{j=1}^N x_{ij} = \alpha$ holds (with $x_{ij}$ depending on $\lambda$). Such a $\lambda$ can be determined by using the binary search---see \cite{bertsimas}  for the example of its use.

Procedure is repeated sequentially for all rows until $X$ has stable form (small differences between iterations). 
According to the suggestions stated in \cite{bertsimas} we set $\alpha=2$. Non--zero (or slightly above zero) entries of obtained matrix $X$ form base for generating initial patterns of the problem (\ref{eq::WTABBfcelu})--(\ref{eq::WTABBconstraint3}). To obtain the initial patterns, we employed sampling based on matrix $X$ in such a way that the set of initial patterns satisfied the constraints of the problem (\ref{eq::WTABBfcelu})--(\ref{eq::WTABBconstraint3}).

\subsubsection{The pricing subproblem as a convex mathematical programming problem}

A second key element of the approach proposed by Bertsimas et al. \cite{bertsimas} is solving the pricing subproblem using techniques for constrained convex optimization. As derived in Subsection \ref{sec::WTAGGlinearization} the pricing subproblem for the RMP problem (\ref{eq::WTABBfcelu})--(\ref{eq::WTABBconstraint3}) is
{\color{black}
\begin{equation}
    \min_{x\in \{0,1\}^M} \left [ V_j e^{s_j^Tx} -\lambda^Tx - \mu_j\right ]. \label{PPcV2_changed}
\end{equation}
}

{\color{black}
Bertsimas et al. \cite{bertsimas} propose to solve first the relaxation of the problem (\ref{PPcV2}):
\begin{equation}
    \min_{x\in [0,1]^M} \left [ V_j e^{s_j^Tx} -\lambda^Tx - \mu_j\right ]. \label{PPcV2}
\end{equation}
}
and then to transform its continuous solution to a pattern of the problem (\ref{eq::WTABBfcelu})--(\ref{eq::WTABBconstraint3}).

To this end we introduce the Lagrange function for the problem (\ref{PPcV2})
\begin{equation}
\mathcal{L}(x, \alpha, \beta) = V_j e^{s^T x} - \lambda^Tx - \mu_j -\alpha^T x + \beta^T (x-1)\nonumber
\end{equation}
and then the KKT conditions: $0 \leq x_i \leq 1$, $\alpha_i \geq 0$, $\beta_i \geq 0$, $\alpha_i x_i = 0$, $\beta_i (x_i-1) = 0$  and $V_j s_i e^{s_j^T x} - \lambda_i - \alpha_i + \beta_i =0$, $i=1,2,\ldots,M$. 

One can observe (see \cite{bertsimas} for details):
\begin{enumerate}
\item If $e^{s^T x} < \frac{\lambda_i}{V_j s_i}$ then $x_i = 0$
\item If $e^{s^T x} > \frac{\lambda_i}{V_j s_i}$ then $x_i = 1$
\item If $e^{s^T x} = \frac{\lambda_i}{V_j s_i}$ then $x_i = \frac{1}{s_i}\left( \ln\left( \frac{\lambda_i}{V_j s_i}\right) - \sum_{k: \frac{\lambda_k}{V_j s_k} < \frac{\lambda_i}{V_j s_i}} s_k\right)$
\end{enumerate}
So it is sufficient to sort indices $i$ by values $\frac{\lambda_i}{V_j s_i}$, and then iterate over them conditioning $e^{s^T x}$ being equal to $\frac{\lambda_k}{V_j s_k}$. If for given position $k$ : $0< x_k < 1$, iteration terminates $x_k$ is set and all entries of $x$ corresponding to indices lower than $k$ are set to $1$. On the other hand if $x_k > 1$ and $x_{k+1} <0$ holds, procedure terminates with $x_i=1$ for indices lower than $k+1$ and $0$ for the rest---resulting in integer solution for relaxed pricing problem. Since at most one entry of relaxed pricing could be non--integer, obtaining proper solution is possible with simple Branch--and--Bound method---repeating iteration skipping branched indices and adjusting sum which is used in the third case. 

\subsubsection{Adding cuts based on cliques}
\label{sec:cuttingplanes}
New patterns are generated by solving a relaxed problem derived from the problem (\ref{eq::WTABBfcelu})–-(\ref{eq::WTABBconstraint3}). The resulting solution to the relaxed problem contains values of variables $y_j^l$ that are not integers (although, typically, the vast majority of $y_j^l$ variables take on values $0$ or $1$). 

To accelerate the process of obtaining integer values of these variables (e.g., using branch--and--bound techniques), it is advantageous to introduce additional constraints into the problem defined by (\ref{eq::WTABBfcelu})–-(\ref{eq::WTABBconstraint3}); these constraints reduce the search space while guaranteeing that the solution to the original problem lies within that restricted space. 

\begin{figure}[ht]
\centering
\includegraphics[width=0.45\linewidth]{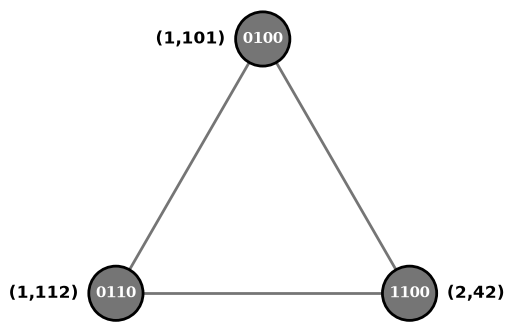} 
\caption{An example of a clique used to generate slashes. The patterns are located in the center of the circles (e.g., 0110 means that pattern no.112 assigned to target no.1 has four weapons that can be used, but in this pattern, only two of them are assigned to the target).} 
\label{fig::clique}
\end{figure}

These inequalities which restrict the search space are called {\it valid inequalities}, among them there {\it clique cuts} which can be applied to the problem (\ref{eq::WTABBfcelu})–(\ref{eq::WTABBconstraint3}). These valid inequalities are derived from the graph perspective of the problem constraints. Suppose that ${\cal G}$ is a graph representing the problem (\ref{eq::WTABBfcelu})-–(\ref{eq::WTABBconstraint3}). Each pattern is represented by a node in the graph ${\cal G}$. If patterns are assigned to the same target $j$, then there exists an edge connecting the nodes representing these patterns. Furthermore, two nodes are connected by an edge if they represent patterns that share a weapon. Suppose that we have a feasible solution to the problem (\ref{eq::WTABBfcelu})-–(\ref{eq::WTABBconstraint3}) which is represented by the graph ${\cal G}$. Then we say that its subset of nodes formulate a clique if any two nodes of the subset are connected---see Fig. \ref{fig::clique}

New cuts are added to the RMP problem only after all subproblems have been examined and it has been determined that no new patterns can be generated. If the RMP solution is not a vector of integers, then additional cuts can be added. 
The cut is always defined on the set $y_{j}^l$, which contains only fractional values ($0 < y_{j}^l < 1$). 
The new constraint has the form:
\begin{equation}
\sum\limits_{(j,l)} {y_{j}^l = 1},\ \forall\ \{(i,j,l):\ 0<y_{j}^l<1 \land b_{ij}^l=1\}\nonumber
\end{equation}

This means that if there are patterns $l$ used to destroy target $j$ for which the value $0<y_j^l<1$ holds, the sum of these $y$ variables must equal one. The operation of adding a cut is marked 'Add Cutting Planes' in Fig. \ref{fig::CGscheme}.

\subsubsection{Mapping RMP relaxation to Integer Assignment}

An effective way to accelerate the convergence of the procedure proposed by Bertsimas et al. \cite{bertsimas} toward a near-optimal solution is to employ the following heuristic: if a small subset of variables does not take on values of $0$ or $1$, these variables are assigned a value of 1 once a specific threshold is exceeded. While this approach does not guarantee an optimal solution, the resulting solution is typically very close to the optimum. In Figure 2, this process is labeled "Update Upper Bound."

The final stage of transforming the solution of the relaxed RMP into an integer solution involves selecting a branching method, which allows for restricting the search space within the relaxed RMP. As always, a well-designed branching rule is essential to avoid analyzing an excessive number of nodes in the Branch-and-Bound tree; naive branching rules can lead to inefficient tree partitioning, rendering the solution process computationally infeasible.

One method used in our work is Most Fractional Variable approach, presented by Land et al. \cite{land1960automatic}, which selects the variable whose fractional part is closest to 0.5. It is very easy and fast to compute, but often generates large search trees. Another tested method was Strong Branching (used to solve the TSP problem by Applegate et al. \cite{applegate1995finding}), which tests several candidate fractional variables. It runs a small number of the simplex method  iterations on each and measures the actual decrease in the objective function value. Based on the calculations, it selects the variable that gives the best improvement. The disadvantage of this method is the significant load on computer processors. However, this method did not yield the expected performance improvement.

Bertsimas et al. \cite{bertsimas} use a branching method, which belongs to the User--defined Priorities group, these methods require expert knowledge of the problem structure. In our method we branch on the linear combinations of the $y_{j}^l$ (see \cite{bertsimas} for details). The procedure requires only calculating two sums for each target and weapon assigned to this target for which the value of $y_{j}^l$ is not an integer. As with the Most Fractional Variable method, a new branch of the tree is created for the variable that is closest to the value of 0.5. In Fig. \ref{fig::CGscheme}, this branching process is marked as "Add New Branch to RMP".

\subsection{Dynamic WTA problem} \label{sec::DynamicWTA}

There are several approaches to dynamic WTA problems (\cite{murphey2000}). Since the approach to dynamic WTA problems proposed in this paper is most similar to the method based on dynamic programming equations, we briefly describe their application in the context of WTA problems. Next, we demonstrate how our method for solving the dynamic WTA problem can be linked to the one utilizing Bellman equations.

In the case of a dynamic WTA problem, the decision-making process regarding the weapon assignments unfolds over time. To simplify the analysis, we assume that this process takes place over $T$ time intervals, each sufficiently long to allow for the employment of weapon systems and the assessment of their effectiveness.
We therefore assume the existence of an information feedback loop, which should ensure more effective use of combat assets, as enemy assets that have already been destroyed will not be targeted in subsequent time intervals. Furthermore, we can expect that threat monitoring will enable the updating of WTA problem parameters, resulting in more effective threat engagement. This strategy is known as "shoot--look--shoot."

The shoot-–look-–shoot strategy involves selecting a subset of weapons in the first time interval and assigning them to targets based on probability estimates available at that interval. Subsequently, a damage assessment is performed, and in the second time interval, targets that have not been neutralized are considered for destruction using the remaining weapons. This procedure is repeated over the subsequent $T-2$ time intervals.

To introduce the dynamic programming equations that form the basis for solving the dynamic WTA problem, we must introduce the following notation: $T$--the number of time intervals; $p_{ij}(t)$--the probability of the $i$th target desctruction by the $j$th weapon at time interval $t$, $i=1,2,\ldots,N$, $j=1,2,\ldots,M$, $t=1,2,\ldots,T$. By $q_{ij}(t)=1-p_{ij}(t)$ we denote the corresponding survival probability. 

For the description of the the Belmann's equations we also need decision variables which can be introduced without the reference to time interval index, the whole solution will include that index:
\begin{eqnarray}
    &{\displaystyle x_{ij}(t) = \left \{ \begin{array}{ll}
    1 & {\rm if}\ {\rm weapon}\ i\ {\rm is}\ {\rm assigned}\ {\rm to}\ {\rm target}\ j\ {\rm on}\ {\rm interval}\ t\nonumber \\
    0 & {\rm otherwise} \nonumber 
    \end{array} \right .}\nonumber \\
    &{\displaystyle t=1,2,\ldots,T}\label{decVariables}
\end{eqnarray}
By $X(t)\in {\cal R}^{M\times N}$ we denote the engagement matrix on the interval $t$, $t=1,2,\ldots,T$

In order to describe the Bellman's equations it is necessary to introduce the concept of the state of the $T$ stages dynamical system. It consists of two vectors: the first indicates {\it state of targets}--$a(t)\in {\cal R}^N$, $t=1,2,\ldots,T$ the second one {\it state of weapons}--$w(t)\in {\cal R}^M$, $t=1,2,\ldots,T$:
\begin{eqnarray}
    &{\displaystyle a_{j}(t) = \left \{ \begin{array}{ll}
    1 & {\rm if}\ {\rm target}\ j\ {\rm survived}\ {\rm on}\ {\rm intervals}\ 1,2,\ldots,t-1 \nonumber \\
    0 & {\rm otherwise} \nonumber 
    \end{array} \right . , }\nonumber \\
    &{\displaystyle j=1,2,\ldots,N}\label{targetState}
\end{eqnarray}
\begin{eqnarray}
    &{\displaystyle w_{i}(t) = \left \{ \begin{array}{ll}
    1 & {\rm if}\ {\rm weapon}\ i\ {\rm was\ not}\ {\rm used}\ {\rm on}\ {\rm intervals}\ 1,2,\ldots,t-1 \nonumber \\
    0 & {\rm otherwise}  \nonumber  
    \end{array} \right . ,}\nonumber \\
    &{\displaystyle i=1,2,\ldots,M} \label{weaponState}
\end{eqnarray}

Notice that $a(t)$ is a random variable since for the purpose of deriving the decision rule $(a(t)^T,w(t)^T)^T = s(t) \rightarrow X(t)$ we consider all possible values of the state vector $s(t+1)$ and they can be described by the conditional distribution 
\begin{eqnarray}
   &{\displaystyle P(a_j(t+1) = k\  \mid s(t)) = k \prod_{i=1}^M \left (1-p_{ij} (t)\right )^{x_{ij}(t)} + (1-k)\left \{ 1 -  \prod_{i=1}^M \left (1-p_{ij} (t)\right )^{x_{ij}(t)} \right \} } \nonumber \\
   &{\displaystyle k=0,1,\ j=1,2,\ldots,N,\ t=1,2,\ldots T} \label{aDistribution}
\end{eqnarray}

On the other hand weapon part of the state vector is a deterministic vector evaluated according to the equation
\begin{eqnarray}
    &{\displaystyle w_i(t) = 1 - \sum_{l=1}^{t-1} \sum_{j=1}^N x_{ij}(l),\ i=1,\ldots,M.}\label{wDistribution}
\end{eqnarray}
Since vectors $w$ are deterministic we can introduce the notation $I(w)= \{i=1,2,\ldots,M:\ w_i =0\}$ which will be employed in the formulation of the Bellman's equations.

Eventually the $t$th stage of the dynamic programming equations look as follow:
\begin{eqnarray}
    &{\displaystyle Q(s(t)) =  \min_{X(t)\in {\cal X}(w(t))} \left [ \sum_{j=1}^N V_j \prod_{i=1}^M \left ( 1-p_{ij}(t)\right )^{x_{ij}(t)} + {\cal E}_{s(t+1)} \left [ Q(s(t+1))\mid s(t), X(t)\right ]  \right ],} \label{eqBellman}
\end{eqnarray}
where ${\cal E}_{u}[Y(u)\mid z]$ is the conditional expectation of $Y[u]$, and
\begin{eqnarray}
    &{\displaystyle {\cal X}(w) = \left \{ X\in {\cal R}^{M\times N}:\ \sum_{j=1}^N x_{ij} = 1,\ i=1,2,\ldots,M,\ x_{ij}\in \{0,1\},\ i\notin I(w),\ x_{ij}=0,\ i\in I(w),\ j=1,2,\ldots,N \right \}.}\label{Xconst}
\end{eqnarray}

The Bellman's equations enable us to establish the decision rule which states the feedback law connecting the current state measurement $s(t)$  with the engagement matrix $\bar{X}(t)$ which is the solution of the Bellman's equations:
\begin{eqnarray}
    &{\displaystyle \bar{X}(s(t)) = f(s(t)) = {\rm arg}\min_{X(t)\in {\cal X}(w)} \left [ \sum_{j=1}^N V_j \prod_{i=1}^M \left ( 1-p_{ij}(t)\right )^{x_{ij}(t)} + {\cal E}_{s(t+1)} \left [ Q(s(t+1))\mid s(t), X(t)\right ]  \right ],} \label{lawFeedback}
\end{eqnarray}

The cost of constructing the feedback function $f$ is enormous because it involves examining the entire state space of the variable $s$; considering only the possible values of the variable $a$, this means the cost is proportional to $2^N$.

However, using the Bellman's equations we not only create the feedback law but also find a global solution to the dynamic WTA problem. Since we know initial conditions for the discrete time dynamical system ($s(1)$) we can construct the strategy of using available weapons at each considered interval---$\bar{X}(1) = f(s(1))$, next we assess the effects of employing weapons due to assets $\bar{X}(1)$, obtaining the current system state $s_c(2)$, and evaluate $\bar{X}(2) = f(s_c(2))$. Proceeding in this manner, at time $t$ we obtain $\bar{X}(t) = f(s_c(t))$, where $s_c(t)$ represents the assessment of damage and the availability of fire assets at time $t$. In this way we create shoot--look--shoot strategy for all intervals $t=1,2,\ldots,T$.

Notable feature of the approach to the dynamic WTA problem is that the implemented weapon employment strategy utilizes actual observations of the system state, whereas the determination of the feedback rule is based on the theoretical probability distribution of the system state. As a result, this approach may lead to situations in which weapons being used against targets that no longer exist.

Due to the high computational cost associated with applying dynamic programming equations to solve the WTA problem, the literature lacks studies presenting computational results for this approach at problem sizes that are practical—for instance, $N \geq 100$ and $M \geq 200$ (\cite{Ahner2015}, \cite{Li2006}). Instead, the Bellman equation methodology is used to derive qualitative results regarding the dynamic approach to the WTA problem—for instance, by demonstrating, through simplified problem scenarios, the superiority of dynamic strategies over strategies based on solving static WTA problems. These simplified cases involve, among other assumptions, setting the destruction probabilities to $p_{ij}(t)=p(t)$ (\cite{Ahner2015}, \cite{murphey2000}, \cite{Li2006}), see also other approaches to the dynamic WTA problem that make similar simplifying assumptions regarding the values of the parameters $p_{ij}$— \cite{Glazebrook2004}, \cite{Burr1985}, \cite{Sikanen2008}, \cite{Soland1987}.

Dynamic programming equations are also used to derive computational methods that are significantly more efficient in implementing weapon engagement strategies—strategies that can be viewed as resulting from the approximate solution of Bellman equations (\cite{Bertsekas2000}, \cite{Davis2017}). In \cite{Bertsekas2000} optimal long term cost $Q(s(t+1))$ is approximated by the function $\tilde{Q}(s(t+1),r)$ where $r$ is a vector of parameters. In \cite{Bertsekas2000} two architectures of a parametric class of functions are considered.

In the case of the first architecture, $\tilde{Q}(s, r)$ is represented as a neural network where the inputs are $s$ and $r$ are its parameters. To determine network parameters, techniques are employed that involve generating solutions to Bellman equations based on approximations of the optimal values $\tilde{Q}(s,r)$, followed by adjustments to the vector $r$.

In the case of the second architecture, called {\it feature extraction mapping}, the state of the system $s$ is transformed into a feature vector
\begin{eqnarray}
    &{\displaystyle g(s) = \left (g_1(s),g_2(s),\ldots,g_p(s)\right ), } \nonumber
\end{eqnarray}
which is then used to represent the function approximating the optimal values in accordance with the equation
\begin{eqnarray}
    &{\displaystyle \tilde{Q}(s,r) = r_0 + \sum_{i=1}^p r_i g_i(s).}\nonumber 
\end{eqnarray}
The values of the parameter vector $r$ are subsequently determined during the training process.

The dynamic WTA problem addressed in \cite{Bertsekas2000} is formulated in a much more general way than the version we present in this article. In particular, the defender must also decide how many weapon to employ against the current attack and how many to reserve for future attacks. The authors present methods from the class of reinforcement learning algorithms known as neuro-dynamic programming—described in \cite{BertsekasT1996}—which are based on approximate policy iteration (API) techniques \cite{Bertsekas2011}. 

In \cite{Davis2017}, another variant of approximate dynamic programming applied to the dynamic WTA problem is considered. The optimal value function is approximated using a linear combination of feature functions, with parameters determined via linear regression methods. However, the approximation of the optimal value function takes place in the post--decision state space. This allows the Bellman equations to be reduced to their deterministic versions (\cite{VanRoy1997}), which significantly lowers the computational cost of solving them.

Methods for solving dynamic WTA problems based on approximate dynamic programming may, with the ongoing development of reinforcement learning techniques, become viable approaches for finding approximate strategies to solve such problems in the future. Simulation results obtained so far for the application of these techniques indicate that dynamic WTA problems with $M \leq 60$ and $N \leq 60$ require several hours of computation time (\cite{Bertsekas2000}).

The method for solving the dynamic WTA problem presented in this paper relies on the MPC approach. This is a dynamic system control method—also known as the receding horizon control—in which the current control action is determined by solving an open-loop optimal control problem over a finite horizon in real time (at each sampling instant), using the system's current state as the initial state (\cite{Mayne2000}, \cite{Clarke1994}, \cite{Keerthi1988}).

What distinguishes the MPC approach from other control methods is the necessity to solve a constrained optimal control problem over a finite horizon—based on the system's current state—in real-time ({\it on-line}) at each time interval $t$. Instead of determining a feedback strategy {\it off-line} (which ensures optimal control for all states)—a process requiring the solution of a Hamilton--Jacobi--Bellman differential or difference equation (associated with dynamic programming) and representing a significantly more difficult task—one solves an open--loop optimal control problem in real--time, using the controlled system's current state as the initial state; this constitutes a mathematical programming problem (specifically, an integer optimization problem in our case). From this perspective, the MPC method differs from other control methods solely in its implementation, which requires that the open--loop optimal control problem be solvable within a reasonable timeframe (relative to the system's dynamics).

When the MPC approach is applied to the dynamic WTA problem the following integer optimization problem has to be solved {\it on-line} on interval $t$
\begin{eqnarray}
    &{\displaystyle  \min_{X \in {\cal X}_t}\sum_{l=t}^T \sum_{j\in {\cal A}_t}V_j \prod_{i\in {\cal W}_t}   \left (1-p_{ij}(l)\right )^{x_{ij}(l)}   }\label{MPCa}
\end{eqnarray}
subject to the constraints
\begin{eqnarray}
    &{\displaystyle    \sum_{l=t}^T\sum_{j\in {\cal A}_t} x_{ij}(l) = 1,\ i\in {\cal W}_t,}\label{MPCb} \\
    &{\displaystyle x_{ij}(l)\in \{0,1\},\ j\in {\cal A}_t,\ i\in {\cal W}_t,\ l=t,t+1,\ldots,T ,}\label{MPCd}
\end{eqnarray}
where
\begin{eqnarray}
    {\cal X}_t & = & \left  \{x_{ij}(l),\ j\in {\cal A}_t,\ i\in {\cal W}_t,\ l=t,t+1,\ldots,T \right\}  \label{MPCna} \\
    {\cal A}_t & = & \left \{j=1,2,\ldots,N:\ a_j(t) = 1\right \} \label{MPCnb} \\
    {\cal W}_t & = & \left \{i=1,2,\ldots,M:\ w_i(t) = 1\right \}. \label{MPCnc}
\end{eqnarray}
Here, ${\cal A}_t$ and ${\cal W}_t$ are {\it on-line} measurements of the state $s$ of the system.

The adopted dynamic model is deterministic and takes the following form
\begin{eqnarray}
    a_j(t+k) & = & a_j(t),\ j=1,2,\ldots,N,\
    k=1,2,\ldots,T-t \label{MPCsa} \\
    w_i(t+k) & = & w_i(t) - \sum_{l=1}^k\sum_{j\in {\cal A}_t} x_{ij}(l),\ i=1,2,\ldots,M, \  k=1,2,\ldots,T-t. \label{MPCsb}
\end{eqnarray}
Notice that $a_j(t+k)=0$, $j\notin {\cal A}_t$, $k=0,1,\ldots T-t$, also $w_i(t+k)=0$, $i\notin {\cal W}_t$, $k=0,1,\ldots T-t$.

The actual evaluation of state $a$ is as follows,
\begin{eqnarray}
    &{\displaystyle       a_j(t+k)  =  a_j(t) - \sum_{l=1}^k{\rm ind}\left (j\in {\cal A}_l\right ),\ j=1,2,\ldots,N,\
    k=1,2,\ldots,T-t    }  \label{MPCsaa}
\end{eqnarray}
Here, 
\begin{eqnarray}
    &{\displaystyle {\rm ind}\left (j\in {\cal A}_l\right )  = \left \{ \begin{array}{ll}
    1 & {\rm if}\ j\in {\cal A}_l \nonumber \\
    0 & {\rm otherwise} \nonumber \end{array} \right .}
\end{eqnarray}

Suppose that $\bar{X}$ is a solution of the problem (\ref{MPCa})--(\ref{MPCd}), it means that
\begin{eqnarray}
    &{\displaystyle  \bar{X} = \arg\min_{X \in {\cal X}_t}\sum_{l=t}^T \sum_{j\in {\cal A}_t}V_j \prod_{i\in {\cal W}_t}   \left (1-p_{ij}(l)\right )^{x_{ij}(l)},   }\label{MPCArg}
\end{eqnarray}
then $\bar{X} = \left ( \{\bar{x}_{ij}(t)\}_{j\in {\cal A}_t,i\in {\cal W}_t} , \{\bar{x}_{ij}(t+1)\}_{j\in {\cal A}_t,i\in {\cal W}_t},\ldots,\{\bar{x}_{ij}(T)\}_{j\in {\cal A}_t,i\in {\cal W}_t}\right )$ and only the first component $\{\bar{x}_{ij}(t)\}_{j\in {\cal A}_t,i\in {\cal W}_t}$ is applied, for the purpose of the subsequent presentation we assign $\{x_{ij}^w(t)\}_{j\in {\cal A}_t,i\in {\cal W}_t}=\{\bar{x}_{ij}(t)\}_{j\in {\cal A}_t,i\in {\cal W}_t}$

Next, at the following interval $t+1$, the result of $\{\bar{x}_{ij}(t)\}_{j\in {\cal A}_t,i\in {\cal W}_t}$ engagement is observed to define ${\cal A}_{t+1}$ and ${\cal W}_{t+1}$ for the optimization problem to follow. We repeat this procedure up to the interval $T$. {\it Algorithm 1} describes the MPC approach to solving the dynamic WTA problem.

\begin{algorithm}[H]
\caption{MPC solution to dynamic WTA}
\KwIn{$N$,\, $M$,\, $\{V_j\}_{j=1}^N$,\, $p_{ij}(t)$,\ $j=1,\ldots,N$,\ $i=1,2,\ldots,M$,\ $t=1,2,\ldots,T$,\ $a_j=1$,\ $j=1,2,\ldots,N$,\ $w_i=1$,\ $i=1,2,\ldots,M$ }
\KwOut{ $ \{x_{ij}^w(t)\}_{j\in {\cal A}_t,i\in {\cal W}_t} $,\ $t=1,2,\ldots,T$}   
\For {$t \in \{1..T\}$ (moving horizon)}{
Evaluate sets ${\cal A}_t$,\ ${\cal W}_t$\;
Solve the problem (\ref{MPCa})--(\ref{MPCd})\;
Let $\bar{X} = \left ( \{\bar{x}_{ij}(t)\}_{j\in {\cal A}_t,i\in {\cal W}_t} , \{\bar{x}_{ij}(t+1)\}_{j\in {\cal A}_t,i\in {\cal W}_t},\ldots,\{\bar{x}_{ij}(T)\}_{j\in {\cal A}_t,i\in {\cal W}_t}\right )$ be the solution\;
Substitute $\{\bar{x}_{ij}(t)\}_{j\in {\cal A}_t,i\in {\cal W}_t}$ for $\{x_{ij}^w(t)\}_{j\in {\cal A}_t,i\in {\cal W}_t}$;
}
\Return $X^w = \left ( \{x^w_{ij}(t)\}_{j\in {\cal A}_t,i\in {\cal W}_t} , \{x^w_{ij}(t+1)\}_{j\in {\cal A}_{t+1},i\in {\cal W}_{t+1}},\ldots,\{x^w_{ij}(T)\}_{j\in {\cal A}_T,i\in {\cal W}_T}\right )$;

\end{algorithm}

The WTA operational method, employing a moving horizon strategy, requires solving a constrained integer optimization problem at each interval. If a given weapon $i$ is used to the target $j$ on the interval $t$, it cannot be targeted at other intervals from the set $[t+1,T]$. If target $j$ has been aimed during interval $t$, whether it has been destroyed will be verified during the observation made in interval $t+1$, and index $j$ will not be included in the set ${\cal A}_{t+1}$ if the attack on it has been successful. This is the reverse of the situation encountered when using dynamic programming equations to determine target neutralization strategies—these equations are based on a Markov decision process that employs target destruction distributions.

\subsection{Linear formulation of the dynamic WTA}
Consider the dynamic WTA problem (formulated for the interval $t=1$)
\begin{equation}
    \min_{x_{ij}(t)} \sum_{t=1}^T \sum_{j=1}^N V_j \prod_{i=1}^M (1-p_{ij})^{x_{ij}(t)} \label{DynamicPata}
\end{equation}
subject to
\begin{equation}
    \sum_{t=1}^T\sum_{j=1}^N x_{ij}(t) = 1,\ i=1,\ldots,M, \label{dynamicPatb}
\end{equation}
\begin{equation}
    x_{ij}(t) \in \{0,1\},\ t=1,\ldots,T \label{dynamicPatc}
\end{equation}

In order to apply Bertsimas approach we have to introduce patterns which will be used in a column generation method. In Bertsimas et al. \cite{bertsimas} for each $j=1,\ldots,N$ we define a vector which indicates which weapons are assigned to the target $j$. Here, the assignment of weapons is for the pair $(j,t)$, $j=1,\ldots,N$, $t=1,\ldots,T$. 

We also define values 
\begin{equation}
    q^l_j(t) = V_j \prod_{i=1}^M (1-p_{ij}(t))^{b^l_{ij}(t)}\label{dynamicPatd}
\end{equation}
then the RMP problem can be stated as
\begin{equation}
    \min_{y_j^l(t)} \sum_{t=1}^T \sum_{j=1}^N \sum_{l=1}^L q^l_j(t) y^l_j(t)\label{dynamicPate}
\end{equation}
subject to
\begin{equation}
\sum_{t=1}^T\sum\limits_{j = 1}^N \sum\limits_{l = 1}^L b_{ij}^l(t) \cdot y_{j}^l(t) = 1,\ \forall i = 1,\ldots,M
\label{dynamicPatf}
\end{equation} 
\begin{equation}
\sum\limits_{l = 1}^L y_{j}^l(t) = 1,\ \forall j = 1,\dots,N,\ t=1,\ldots,T
\label{dynamicPatg}
\end{equation} 
\begin{equation}
y_{j}^l(t) \in \{0,1\},\ \forall j = 1,\ldots,N,\ t=1,\ldots,T,\ l = 1,\ldots,L.
\label{dynamicPath}
\end{equation} 

Suppose that $\mu(t) \in {\cal R}^N$ are the Lagrange multipliers of the dual problem associated with the constraints (\ref{dynamicPatc}) and $\lambda \in {\cal R}^M$ are the Lagrange multipliers associated with the constraints (\ref{dynamicPatb}). According to Bertsimas et al. \cite{bertsimas} the reduced cost of the pattern $l$ evaluated for the target $j$ is expressed as
\begin{equation}
    q_j^l(t) - \mu_j(t) - \sum_{i=1}^M b^l_{ij}(t) \lambda_i \label{dynPath}
\end{equation}
If the reduced cost, for a given pattern, is negative then incorporating this pattern to the current set of patterns can reduce value of the objective function of the relaxed RMP. Therefore, for each target $j$, $j=1,\ldots,N$, $t=1,\ldots,T$ we solve the pricing subproblem
\begin{equation}
    \min_{x_{ij}\in \{0,1\}^M} \left [ V_j \prod_{i=1}^M (1-p_{ij}(t))^{x_{ij}} -\lambda^Tx - \mu_j(t)\right ]. \label{dynPatu}
\end{equation}

The subproblems indexed by $(j,t)$, $j=1,\ldots<N$, $t=1,\ldots,T$ can be solved by the Bertsmas et al. \cite{bertsimas} procedure proposed for their Pricing Problem.

\section{Results}
\label{sec:Results}
Testing environment specifications:
\begin{itemize}
    \item Intel(R) i7-9700(R) CPU @ 3.0 GHz (Cores: 8 Threads: 8)
    \item 32.0 GB RAM
    \item 512 GB SSD storage
    \item 64-bit Windows 10 operating system.
\end{itemize}

Regarding the processor used for the calculations, experiments were conducted utilizing the available number of physical cores (4–6). The tables present the calculation results for six cores. The number of physical cores employed affects the time required to solve the subproblems, whereas it has a lesser impact on the time needed to solve the RMP task.

The algorithms were tested on three test sets available online. The first set was developed by Ahuja (test cases prepared by Lou and Chen - test cases obtained from the authors), the second set was developed by Sonuc (test cases: \url{https://data.mendeley.com/datasets/jt2ppwr62p/2}), and the third set was developed by Toffolo and Pavlikov (GitHub: \url{https://github.com/tuliotoffolo/wta}, see \cite{andersen}). Each of these sets contains test cases of varying sizes, including the most computationally difficult cases for WTA algorithms. These are cases where the number of destructive agents $M$ is greater than the number of threats being countered $N$.

\subsection{Results for static WTA problem}
\label{sec::WTAstatic}

The experiments were conducted using two possible approaches. The first approach did not use a dedicated initialization algorithm, which was developed based on the guidelines outlined by Bertsimas et al. \cite{bertsimas}. In the second approach, the initialization algorithm used in each case was developed to generate initial solutions.

The attached results demonstrate significant speedups when using the initialization algorithm based on Large Neighborhood Search techniques used by Bertsimas.

It's important to note one important aspect that accelerates problem solving, particularly when searching for an optimal solution. The basic version of the WTA algorithm assumes the use of all available resources to destroy threats. This is not implemented in practice, especially when the number of resources to be destroyed significantly exceeds the number of incoming threats. No one fires everything they have for economic reasons. In such a case, the speedup is most noticeable for the algorithm based on the subproblem, where the subproblem is an integer problem, where the objective function is linearized according to Andersen's guidelines.

\begin{table}[!h]
  \caption{Results for instances prepared by Lu\&Chen. Simple
initialization algorithm included.}
   \centering
   \resizebox{\textwidth}{!}{%
\begin{tabular}{lrrrrrrrrrrrr}
\toprule
  & \multicolumn{2}{r}{time [s]} & \multicolumn{2}{r}{score} &
\multicolumn{2}{r}{patterns} & \multicolumn{2}{r}{RMP iter} \\
  & linear WTA & LuChen & linear WTA & LuChen & linear WTA & LuChen &
linear WTA & LuChen \\
  & sec.\ref{sec::WTAGGlinearization} & sec.\ref{sec::WTAGGChen} &  &  &
&  &  & \\
instance &  &  &  &  &  &  &  &  & \\
\midrule

Ahuja\_5\_5 & 0.007 & 0.003 & 87.557 & 43.139 & 10 & 26 & 1 & 4 \\
Ahuja\_10\_10 & 0.038 & 0.004 & 94.347 & 94.347 & 82 & 54 & 9 & 5 \\
Ahuja\_10\_5 & 0.038 & 0.012 & 7.862 & 7.764 & 58 & 82 & 15 & 17 \\
Ahuja\_10\_20 & 0.043 & 0.003 & 538.107 & 538.107 & 112 & 73 & 6 & 4 \\
Ahuja\_20\_80 & 0.141 & 0.005 & 3508.970 & 3508.970 & 443 & 296 & 8 & 5
\\
Ahuja\_20\_40 & 0.139 & 0.003 & 1126.500 & 1125.950 & 411 & 158 & 12 & 4
\\
Ahuja\_20\_20 & 0.094 & 0.005 & 163.664 & 162.809 & 236 & 109 & 13 & 7
\\
Ahuja\_20\_10 & 0.103 & 0.020 & 12.325 & 12.207 & 221 & 213 & 25 & 25 \\
Ahuja\_40\_40 & 0.237 & 0.007 & 249.929 & 249.929 & 885 & 236 & 24 & 7
\\
Ahuja\_40\_20 & 0.221 & 0.024 & 18.538 & 18.294 & 588 & 392 & 32 & 24 \\
Ahuja\_40\_120 & 0.280 & 0.015 & 4520.830 & 4520.060 & 1164 & 551 & 13 &
8 \\
Ahuja\_40\_80 & 0.229 & 0.010 & 2182.470 & 2182.180 & 907 & 458 & 14 & 7
\\
Ahuja\_80\_320 & 0.719 & 0.082 & 13605.300 & 13605.300 & 2711 & 1484 &
15 & 8 \\
Ahuja\_80\_160 & 0.560 & 0.048 & 4331.090 & 4330.820 & 2361 & 969 & 18 &
9 \\
Ahuja\_80\_80 & 0.727 & 0.023 & 544.681 & 544.681 & 2975 & 444 & 39 & 7
\\
Ahuja\_80\_40 & 0.525 & 0.074 & 30.176 & 29.915 & 1455 & 801 & 41 & 34
\\
Ahuja\_100\_100 & 0.953 & 0.054 & 614.366 & 614.260 & 3969 & 569 & 41 &
10 \\
Ahuja\_100\_200 & 0.868 & 0.075 & 4909.480 & 4909.080 & 3628 & 1106 & 22
& 8 \\
Ahuja\_100\_50 & 0.771 & 0.126 & 37.690 & 37.561 & 2034 & 1028 & 47 & 35
\\
Ahuja\_200\_400 & 2.784 & 0.479 & 10232.800 & 10232.200 & 11762 & 2323 &
34 & 8 \\
Ahuja\_200\_100 & 2.598 & 0.735 & 66.300 & 66.096 & 4951 & 2260 & 62 &
34 \\
Ahuja\_200\_200 & 3.775 & 0.303 & 1280.480 & 1279.900 & 13850 & 1215 &
72 & 9 \\
Ahuja\_400\_400 & 25.384 & 2.250 & 2473.210 & 2473.160 & 49492 & 2305 &
127 & 9 \\
\bottomrule
\end{tabular}%
}
\label{tab::LuChenwithoutInit}
\end{table}

The results presented in Table \ref{tab::LuChenwithoutInit} concern the performance of WTA algorithms on the Ahuja test sets. In this case, no advanced method was used to generate the initial set of patterns (so the Bertsimas method is not presented). The absence of an initial set of patterns—which would define acceptable solutions close to the optimum—means that the individual algorithms require more time to arrive at optimal solutions. For an algorithm that employs an integer programming problem as a subproblem (linear WTA), it is advisable to limit the maximum number of weapons that can be assigned to a given target (to a value not exceeding 4; this approach is also used by Andersen {\cite{andersen}}). This modification does not affect the quality of the solutions obtained. Any minor differences in the objective function values between the various models are a consequence of linearizing the exponential function within the linear programming formulation.
For every problem where a specialized mechanism for generating exact patterns was not employed, the use of additional methods within the column generation procedure resulted in a significant reduction in computation time. This applies in particular to the "Update Upper Bound" step (see \ref{fig::CGscheme}).

\begin{table}[!h]

\caption{Results for instances prepared by Lu\&Chen. Modified initialization algorithm included.}

  \centering

  \resizebox{\textwidth}{!}{%

\begin{tabular}{lrrrrrrrrrrrrr}

\toprule

& \multicolumn{4}{r}{time} & \multicolumn{3}{r}{score} & \multicolumn{3}{r}{patterns} & \multicolumn{3}{r}{RMP iter} \\

& init & linear WTA & LuChen & Bertsimas & linear WTA & LuChen &  Bertsimas & linear WTA & LuChen & Bertsimas & linear WTA & LuChen & Bertsimas\\

&  & sec.\ref{sec::WTAGGlinearization} & sec.\ref{sec::WTAGGChen} & sec.\ref{sec::WTAGGBertsimas}  &  &  &  &  &  &  &  &\\

instance &  &  &  &  &  &  &  &  &  &  &  &  &\\

\midrule

Ahuja\_5\_5 & 0.000 & 0.007 & 0.009 & 0.003 & 43.139 & 43.139 & 43.139 & 27 & 25 & 1 & 1 & 6 & 2 \\

Ahuja\_10\_10 & 0.000 & 0.012 & 0.005 & 0.004 & 94.347 & 94.347 & 94.347 & 73 & 66 & 5 & 2 & 6 & 3 \\

Ahuja\_10\_5 & 0.000 & 0.016 & 0.014 & 0.004 & 7.764 & 7.764 & 7.764 & 162 & 75 & 18 & 3 & 19 & 14\\

Ahuja\_10\_20 & 0.000 & 0.025 & 0.006 & 0.004 & 538.107 & 538.107 & 538.107 & 83 & 98 & 6 & 3 & 8 & 3\\

Ahuja\_20\_80 & 0.001 & 0.081 & 0.008 & 0.005 & 3508.970 & 3508.970 & 3508.966 & 231 & 274 & 47 & 4 & 9 & 3\\

Ahuja\_20\_40 & 0.001 & 0.040 & 0.007 & 0.005 & 1125.950 & 1125.950 & 1125.955 & 177 & 210 & 20 & 3 & 8 & 3\\

Ahuja\_20\_20 & 0.000 & 0.020 & 0.011 & 0.004 & 162.809 & 162.809 & 162.809 & 180 & 147 & 6 & 2 & 14 & 2\\

Ahuja\_20\_10 & 0.000 & 0.023 & 0.018 & 0.005 & 12.207 & 12.207 & 12.207 & 305 & 178 & 38 & 4 & 22 & 4\\

Ahuja\_40\_40 & 0.001 & 0.051 & 0.010 & 0.008 & 249.929 & 249.929 & 249.929 & 422 & 310 & 362 & 3 & 9 & 21\\

Ahuja\_40\_20 & 0.001 & 0.027 & 0.023 & 0.009 & 18.508 & 18.294 & 18.294 & 626 & 361 & 96 & 2 & 21 & 6\\

Ahuja\_40\_120 & 0.003 & 0.128 & 0.022& 0.010 & 4520.060 & 4520.060 & 4520.062 & 431 & 635 & 109 & 4 & 12 & 3\\

Ahuja\_40\_80 & 0.002 & 0.097 & 0.016 & 0.009 & 2182.180 & 2182.180 & 2182.177 & 372 & 520 & 55 & 4 & 12 & 4\\

Ahuja\_80\_320 & 0.014 & 0.338 & 0.155 & 0.030 & 13605.300 & 13605.300 & 13605.302 & 1114 & 1693 & 381 & 7 & 16 & 4\\

Ahuja\_80\_160 & 0.007 & 0.196 & 0.088 & 0.021 & 4330.900 & 4330.820 & 4330.818 & 899 & 1627 & 259 & 5 & 16 & 4\\

Ahuja\_80\_80 & 0.006 & 0.086 & 0.047 & 0.017 & 544.722 & 544.681 & 544.681 & 817 & 792 & 70 & 3 & 15 & 3\\

Ahuja\_80\_40 & 0.004 & 0.100 & 0.083 & 0.023 & 29.949 & 29.915 & 29.915 & 1508 & 997 & 377 & 5 & 36 & 7\\

Ahuja\_100\_100 & 0.009 & 0.129 & 0.094 & 0.023 & 614.318 & 614.260 & 614.260 & 997 & 1182 & 53 & 4 & 18 & 4\\

Ahuja\_100\_200 & 0.014 & 0.189 & 0.165 & 0.027 & 4909.250 & 4909.080 & 4909.085 & 1256 & 2072 & 148 & 4 & 18 & 3\\

Ahuja\_100\_50 & 0.006 & 0.142 & 0.136 & 0.031 & 37.581 & 37.561 & 37.561 & 1862 & 1091 & 356 & 6 & 30 & 7\\

Ahuja\_200\_400 & 0.055 & 0.438 & 2.603 & 0.092 & 10232.400 & 10232.200 & 10232.154 & 2836 & 11904 & 282 & 5 & 38 & 5\\

Ahuja\_200\_100 & 0.024 & 0.356 & 1.088 & 0.147 & 66.198 & 66.096 & 66.096 & 4298 & 2918 & 859 & 7 & 43 & 7\\

Ahuja\_200\_200 & 0.032 & 0.250 & 0.602 & 0.105 & 1280.160 & 1279.900 & 1279.899 & 2568 & 2512 & 164 & 4 & 18 & 5\\

Ahuja\_400\_400 & 0.131 & 0.868 & 10.310 & 0.539 & 2473.260 & 2473.160 & 2473.155 & 6484 & 13638 & 418 & 5 & 42 & 4\\

\bottomrule

\label{tab::LuChenwithInit}
\end{tabular}
}
\end{table}

The results presented in Table \ref{tab::LuChenwithInit} concern the efficiency of WTA algorithms on Ahuja test sets. In this instance, an advanced method for generating the input pattern set—based on a modified Bertsimas concept—was employed. A significant reduction in solution search time is evident in virtually every case where $M>N$. This is particularly noticeable for the method in which the subproblem is formulated as an integer programming problem. The data presented in the table indicate that the method based on the linearized WTA formulation does not slow down as the dimensionality increases, which is a key conclusion drawn from the study of the method in this scope.
The results presented here are slightly better than those reported by Bertsimas, which is due to the use of newer multi-core i7 processors.

Comparing the results presented in Tables \ref{tab::LuChenwithoutInit} and \ref{tab::LuChenwithInit}, one can observe that the improved initialization algorithm directly affected two elements of the CG algorithm. The first is the number of generated patterns, which is often significantly lower for the algorithm with improved initialization. This frequently translates into fewer iterations required to solve the RMP (see Fig. \ref{fig::CGscheme} – the upper loop encompassing the tasks of adding cutting planes and branching within the *branch-and-bound* method). The reduction in the number of iterations is primarily driven by algorithm steps involving the generation of new clique-based cutting planes (Add Cutting Planes), the updating of the WTA problem's upper bound (Update Upper Bounds), and the modified method for selecting variables to construct the B\&B tree (Add New Branch to RMP).

At the same time, it should be noted that the method for generating initial patterns proposed by Bertsimas is not well-suited for the column enumeration method developed by Chen. In this case, the results deteriorate significantly.

\begin{table}[!h]
  \caption{Results for instances prepared by Sonuc\&Bayir. Simple initialization algorithm included.}
   \centering
   \resizebox{\textwidth}{!}{%
\begin{tabular}{lrrrrrrrrrrrr}
\toprule
& \multicolumn{2}{r}{time [s]} & \multicolumn{2}{r}{score} &
\multicolumn{2}{r}{patterns} & \multicolumn{2}{r}{RMP iter} \\
  & linear WTA & LuChen & linear WTA & LuChen & linear WTA & LuChen &
linear WTA & LuChen \\
  & sec.\ref{sec::WTAGGlinearization} & sec.\ref{sec::WTAGGChen} &  &  &
&  &  &  \\
instance &  &  &  &  &  &  &  &  & \\
\midrule

Sonuc\_wta5 & 0.007 & 0.004 & 48.438 & 48.364 & 10 & 30 & 1 & 5 \\
Sonuc\_wta10 & 0.052 & 0.005 & 96.312 & 96.312 & 94 & 64 & 12 & 7 \\
Sonuc\_wta20 & 0.127 & 0.005 & 142.107 & 142.107 & 301 & 133 & 18 & 6 \\
Sonuc\_wta30 & 0.168 & 0.004 & 248.525 & 248.029 & 504 & 156 & 18 & 5 \\
Sonuc\_wta40 & 0.255 & 0.008 & 305.502 & 305.502 & 926 & 218 & 26 & 7 \\
Sonuc\_wta50 & 0.304 & 0.011 & 353.454 & 353.010 & 1193 & 305 & 25 & 8
\\
Sonuc\_wta60 & 0.396 & 0.013 & 414.276 & 414.222 & 1713 & 343 & 30 & 7
\\
Sonuc\_wta70 & 0.516 & 0.014 & 495.845 & 495.775 & 1975 & 384 & 31 & 6
\\
Sonuc\_wta80 & 0.573 & 0.030 & 530.794 & 530.574 & 2385 & 458 & 31 & 9
\\
Sonuc\_wta90 & 0.773 & 0.031 & 587.816 & 587.816 & 3265 & 499 & 39 & 8
\\
Sonuc\_wta100 & 0.971 & 0.043 & 692.157 & 691.924 & 3710 & 589 & 40 & 8
\\
Sonuc\_wta200 & 4.366 & 0.257 & 1275.150 & 1274.780 & 15363 & 1131 & 79
& 8 \\
\bottomrule
\end{tabular}
}
   \label{tab::SonucwithoutInit}
\end{table}

The results presented in Table \ref{tab::SonucwithoutInit} concern the performance of WTA algorithms on the Sonuc test sets, where $N=M$. This set of cases is significantly easier for WTA algorithms tasked with assigning each weapon to a target. As with the other test sets, no advanced method for generating the initial set of patterns was employed for the Ahuja-type data. The time required to find the optimal solution is longer than when the initial pattern set is carefully selected.
For the Sonuc test sets, each algorithm finds optimal or near-optimal solutions (any potential error depends solely on the linearization accuracy of the WTA subproblem's objective function). The number of RMP iterations is lower than for the Ahuja test sets. The algorithm employing an integer programming formulation for the CG subproblem generates a number of patterns comparable to that obtained for the Ahuja sets.

\begin{table}[!h]

\caption{Results for instances prepared by Sonuc\&Bayir. Modified initialization algorithm included.}

  \centering

  \resizebox{\textwidth}{!}{%

\begin{tabular}{lrrrrrrrrrrrrr}

\toprule

& \multicolumn{4}{r}{time} & \multicolumn{3}{r}{score} & \multicolumn{3}{r}{patterns} & \multicolumn{3}{r}{RMP iter} \\

& init & linear WTA & LuChen & Bertsimas & linear WTA & LuChen & Bertsimas & linear WTA & LuChen & Bertsimas & linear WTA & LuChen & Bertsimas \\

&  & sec.\ref{sec::WTAGGlinearization} & sec.\ref{sec::WTAGGChen} & sec.\ref{sec::WTAGGBertsimas} &  &  &  &  &  &  &  &\\

instance &  &  &  &  &  &  &  &  &  &  &  &  &\\

\midrule

Sonuc\_wta5 & 0.000 & 0.007 & 0.004 & 0.002 & 48.364 & 48.364 & 48.364 & 28 & 28 & 1 & 1 & 6 & 1 \\

Sonuc\_wta10 & 0.000 & 0.020 & 0.006 & 0.004 & 96.312 & 96.312 & 96.312 & 91 & 74 & 11 & 4 & 8 & 4\\

Sonuc\_wta20 & 0.000 & 0.034 & 0.010 & 0.004 & 142.371 & 142.107 & 142.107 & 146 & 176 & 22 & 4 & 12 & 3 \\

Sonuc\_wta30 & 0.001 & 0.035 & 0.009 & 0.005 & 248.029 & 248.029 & 248.029 & 273 & 233 & 6 & 3 & 10 & 2\\

Sonuc\_wta40 & 0.001 & 0.045 & 0.012 & 0.007 & 305.502 & 305.502 & 305.502 & 392 & 330 & 18 & 3 & 11 & 4\\

Sonuc\_wta50 & 0.002 & 0.080 & 0.017 & 0.010 & 353.202 & 353.010 & 353.010 & 492 & 441 & 85 & 5 & 12 & 3\\

Sonuc\_wta60 & 0.003 & 0.080 & 0.029 & 0.011 & 414.455 & 414.222 & 414.222 & 553 & 725 & 40 & 4 & 16 & 4\\

Sonuc\_wta70 & 0.003 & 0.086 & 0.037 & 0.013 & 495.796 & 495.775 & 495.775 & 662 & 732 & 37 & 3 & 15 & 3\\

Sonuc\_wta80 & 0.006 & 0.107 & 0.048 & 0.018 & 530.574 & 530.574 & 530.574 & 867 & 927 & 54 & 4 & 15 & 4\\

Sonuc\_wta90 & 0.005 & 0.098 & 0.086 & 0.020 & 587.901 & 587.816 & 587.816 & 907 & 971 & 52 & 3 & 16 & 3\\

Sonuc\_wta100 & 0.009 & 0.111 & 0.105 & 0.026 & 692.051 & 691.924 & 691.924 & 1171 & 1175 & 57 & 3 & 18 & 4\\

Sonuc\_wta200 & 0.027 & 0.207 & 0.778 & 0.098 & 1274.780 & 1274.780 & 1274.778 & 2776 & 3122 & 147 & 3 & 23 & 5\\

\bottomrule
\end{tabular}

}
  \label{tab::SonucwithInit}
\end{table}

Table \ref{tab::SonucwithInit} presents the results of the WTA algorithms for the case where initial patterns were generated using an improved Bertsimas algorithm. This approach reduces the solution search time. This is particularly evident in the first model, which involves integer subproblems. In this variant, both the solution search time and the number of RMP iterations are reduced, resulting in a decrease in the number of generated patterns. For the presented test cases, an optimal solution can be obtained in under one second in every instance. Furthermore, linear scalability with respect to problem size can be observed for each of the algorithms.
The column enumeration algorithm is also less efficient when it employs the Bertsimas initialization and test data from the Sonuc datasets.

\begin{table}[!h]
  \caption{Results for instances prepared by Tulio\&Toffolo. Simple
initialization algorithm included.}
   \centering
   \resizebox{\textwidth}{!}{%
\begin{tabular}{lrrrrrrrrrrrr}
\toprule
  & \multicolumn{2}{r}{time [s]} & \multicolumn{2}{r}{score} &
\multicolumn{2}{r}{patterns} & \multicolumn{2}{r}{RMP iter} \\
  & linear WTA & LuChen & linear WTA & LuChen & linear WTA & LuChen &
linear WTA & LuChen \\
  & sec.\ref{sec::WTAGGlinearization} & sec.\ref{sec::WTAGGChen} &  &  &
&  &  & \\
instance &  &  &  &  &  &  &  &  &  \\
\midrule
tt\_wta\_50x100 & 0.290 & 0.013 & 1737.540 & 1734.300 & 1280 & 497 & 15
& 7 \\
tt\_wta\_100x200 & 0.784 & 0.082 & 2864.400 & 2862.300 & 3168 & 1130 &
20 & 9 \\
tt\_wta\_150x300 & 0.885 & 0.247 & 3745.580 & 3740.420 & 3858 & 1759 &
16 & 9 \\
tt\_wta\_200x400 & 1.648 & 0.616 & 4973.460 & 4970.360 & 5936 & 2661 &
20 & 10 \\
tt\_wta\_250x500 & 1.971 & 1.766 & 6694.330 & 6689.830 & 6981 & 4068 &
18 & 14 \\
tt\_wta\_300x600 & 3.562 & 2.678 & 8143.210 & 8136.880 & 12056 & 5028 &
25 & 13 \\
tt\_wta\_350x700 & 3.740 & 4.879 & 9059.380 & 9049.830 & 11940 & 6358 &
22 & 15 \\
tt\_wta\_400x800 & 5.316 & 6.216 & 11496.700 & 11491.000 & 16052 & 6362
& 25 & 13 \\
tt\_wta\_450x900 & 6.764 & 9.663 & 12524.400 & 12521.100 & 17924 & 7188
& 26 & 14 \\
tt\_wta\_500x1000 & 7.963 & 14.079 & 13148.300 & 13141.400 & 20001 &
9778 & 26 & 15 \\
\bottomrule
\end{tabular}
}
   \label{tab::TuliowithoutInit}
\end{table}

\begin{table}

\caption{Results for instances prepared by Tulio\&Toffolo. Modified initialization algorithm included.}

  \centering

  \resizebox{\textwidth}{!}{%

\begin{tabular}{lrrrrrrrrrrrrr}

\toprule

& \multicolumn{4}{r}{time} & \multicolumn{3}{r}{score} & \multicolumn{3}{r}{patterns} & \multicolumn{3}{r}{RMP iter} \\

& init & linear WTA & LuChen & Bertsimas  & linear WTA & LuChen & Bertsimas & linear WTA & LuChen & Bertsimas  & linear WTA & LuChen  & Bertsimas\\

&  & sec.\ref{sec::WTAGGlinearization} & sec.\ref{sec::WTAGGChen} & sec.\ref{sec::WTAGGBertsimas} &  &  &  &  &  &  &  &    \\

instance &  &  &  &  &  &  &  &  &  &  &  &  &  \\

\midrule

tt\_wta\_50x100  & 0.004 & 0.227 & 0.021 & 0.009 & 1741.300 & 1734.300 & 1734.302 & 595 & 606 & 143 & 11 & 11 & 4 \\

tt\_wta\_100x200 & 0.021 & 0.247 & 0.121 & 0.021 & 2863.250 & 2862.300 & 2862.302 & 921 & 1421 & 353 & 6 & 14 & 4\\

tt\_wta\_150x300 & 0.038 & 0.405 & 0.469 & 0.039 & 3740.940 & 3740.420 & 3740.418 & 1448 & 2458 & 529 & 6 & 17 & 5\\

tt\_wta\_200x400 & 0.070 & 0.555 & 1.373 & 0.066 & 4973.240 & 4970.360 & 4970.363 & 2058 & 5372 & 867 & 7 & 21 & 5\\

tt\_wta\_250x500 & 0.120 & 0.785 & 3.616 & 0.103 & 6692.730 & 6689.830 & 6689.831 & 2618 & 8847 & 1317 & 8 & 27 & 6\\

tt\_wta\_300x600 & 0.183 & 0.932 & 5.832 & 0.134 & 8139.670 & 8136.880 & 8136.880 & 3018 & 10215 & 1509 & 7 & 28 & 5\\

tt\_wta\_350x700 & 0.241 & 0.988 & 11.075 & 0.182 &9053.380 & 9049.830 & 9049.826 & 3653 & 15584 & 1776 & 6 & 33 & 5\\

tt\_wta\_400x800 & 0.347 & 1.379 & 14.326 & 0.239 &11491.850 & 11491.000 & 11491.045 & 4068 & 14826 & 1999 & 7 & 29 & 5\\

tt\_wta\_450x900 & 0.274 & 1.909 & 26.954 & 0.298 &12522.100 & 12521.100 & 12521.103 & 4264 & 25579 & 2463 & 8 & 38 & 5\\

tt\_wta\_500x1000 & 0.507 & 1.937 & 31.766 & 0.371 &13143.300 & 13141.400 & 13141.439 & 5234 & 21332 & 2652 & 8 & 33 & 5\\

\bottomrule
\end{tabular}
}
  \label{tab::TuliowithInit}
\end{table}

Tables \ref{tab::TuliowithoutInit} and \ref{tab::TuliowithInit} present the results obtained by Weapon-Target Assignment (WTA) algorithms on test sets prepared by Tulio and Toffolo—a collection previously used in the work of Bertsimas \cite{bertsimas}. The test set includes instances with dimensions of 350x700, 400x800, and larger; these represent the most significant challenge for the proposed algorithms.
The results indicate that instances exceeding 350x700 in size do not pose a substantial challenge for any of the algorithms discussed in this paper. The most effective method is the one based on the approach proposed by Bertsimas. Nevertheless, even the algorithm employing linearization solves all tasks in under 2 seconds. It is also worth noting that these methods exhibit linear scalability, even in scenarios where $M > N$.
The test sets prepared by Tulio Toffolo contain the most complex cases, in which the number of weapons is significantly greater than the number of targets. 

For larger problems, the importance of an effective initialization algorithm becomes clearly apparent, as it can reduce the time required to find a solution by several tens of percent. Furthermore, it should be noted that the algorithm by Lu \cite{chen} does not perform efficiently when using the set of initial patterns determined by the Bertsimas method.

\begin{table}[!h]
  \centering
  \caption{A very large-scale instances.}
  \label{tab:largeInstances}
\begin{tabular}{rrrrrr}
\toprule
\multicolumn{3}{r}{instance} & \multicolumn{3}{r}{Bertsimas Branch-Price-and-Cut}
\\
Targets & Unique weapons & Total weapons & time & patterns & RMP iter \\
\midrule
2500  & 1250 & 1250  & \textcolor{black}{1.66}   & \textcolor{black}{19 528} & \textcolor{black}{6} \\
2500  & 1250 & 2500  & \textcolor{black}{4.51}   & \textcolor{black}{29 356} &  \textcolor{black}{8} \\
2500  & 1250 & 3750  & \textcolor{black}{114.38}  & \textcolor{black}{59 977} & \textcolor{black}{68} \\
5000  & 2500 & 2500  & \textcolor{black}{8.90}   & \textcolor{black}{40 238} &  \textcolor{black}{6} \\
5000  & 2500 & 5000  & \textcolor{black}{21.70}  & \textcolor{black}{68 371} &  \textcolor{black}{7} \\
5000  & 2500 & 7500  & \textcolor{black}{1169.19} & \textcolor{black}{135 000} & \textcolor{black}{99} \\
10000 & 5000 & 5000  & \textcolor{black}{40.69}  & \textcolor{black}{84 357} & \textcolor{black}{6} \\
10000 & 5000 & 10000 & \textcolor{black}{183.87} & \textcolor{black}{151 723} & \textcolor{black}{11} \\
\bottomrule
\end{tabular}
\end{table}

The Branch-Price-and-Cut method was applied to very large synthetic instances generated using the procedure developed by Bertsimas \cite{bertsimas}. The results are presented in Table \ref{tab:largeInstances}. To increase the problem's complexity, the option to define weapon types with identical parameters (i.e., the same target destruction probabilities) was introduced. The results obtained by Bertsimas method significantly outperform all methods described in the article in terms of solution effectiveness.
But it is worth noting that the method based on objective function linearization and solving subproblems as MIP tasks determines optimal values in a timeframe that compares favorably to the methods described in \cite{andersen}. It is also efficient when compared to heuristic algorithms. Compared to the method presented by \cite{bertsimas}, it can be observed that as the problem size increases, the ratio of the computation times of the linear WTA method to those of the Bertsimas-based method does not change significantly in favor of the latter.

\subsection{Results for dynamic WTA problem }
\label{sec::WTAdynamic}

This section presents the results for the DWTA (Dynamic Weapon-Target Assignment) problem in comparison to the static WTA problem, covering three consecutive time intervals during which threats were engaged using available weapon systems.

The following assumptions were made for the model:
\begin{itemize}
    \item The set of targets detected in the first time interval is not expanded; that is, no new threats are detected in subsequent intervals.

    \item Targets detected in the first interval have a specific probability of destruction assigned to each weapon.

    \item Whether a target has been destroyed, or remains available for destruction in a subsequent time interval, is determined by the method described below.

    \item If a target is available for destruction in the next time interval, the probability of its destruction by a weapon may increase, decrease, or remain unchanged.
    
\end{itemize}

Comparative tests were conducted for two WTA models: the static WTA model and the dynamic WTA model with an additive objective function (DWTA-A \eqref{dynamicPate}).

The tests involved comparing the statistical outcomes of weapon-to-target allocations for the SWTA model—applied to a static instance (defined as assigning weapons to targets within a single time window)—and the DWTA model—applied to a dynamic instance (covering three time intervals). The dynamic instances were generated from the static instance using four different methods. In the "equal" variant, the probability matrix $p_k$ is replicated across all three time windows, meaning the engagement probabilities remain constant over time. The "random" variant introduces noise to the "equal" variant by adding values drawn from a uniform distribution $[-0.1, 0.1]$.
The "increasing" and "decreasing" variants are modifications of the "random" variant; in these cases, all engagement probabilities either increase or decrease over time (though some specific values may remain unchanged).

In the SWTA approach, weapons are allocated to targets, and the destruction of each target is subsequently verified (a process known as "kill assessment") by drawing a random number from a $U(0,1)$ distribution. If, for a given target $j$, this number is greater than $\prod_{k : x_{kj} =1} (1-p_{kj})$—which represents the probability of survival under the specific allocation—then the target is considered destroyed. The number and values of the targets that survive are then tallied. This sampling procedure is repeated $N_s = 200$ times.

In the DWTA approach, the kill assessment process is carried out after each time interval—three times in the examples presented. In the Receding Horizon method, only those weapons and targets that have not been expended or destroyed in previous time steps are carried over to the next time interval. At each stage, the counts and values of the remaining targets, as well as the available weapons, are tallied. This sampling process is repeated $N_s$ times. The tables compare the statistics from the numerical experiments conducted (including statistics on the DWTA computation time across successive time intervals).

In the following tables summarizing the experiments conducted for the DWTA scenarios, the test results for successive time intervals are separated by a vertical bar "|". The average number of remaining targets and weapons reflects the state after simulating the destruction resulting from the allocations made within a given time inyrtval. Due to the nature of the SWTA (Static Weapon-Target Assignment) component, all available weapons are always consumed. Similarly, if any target remains by the final time interval, the allocation process will utilize all remaining weapons (in line with the objective of minimizing the expected survival value). Consequently, at the end of every instance—even when the expected number of remaining targets is zero—all available weapons are exhausted. It is worth considering a model modification that would allow for the conservation of ammunition; this could be achieved by incorporating the cost of using a weapon to destroy a target into the objective function.
Table \ref{tab::DWTA1} presents calculations for both static and dynamic WTA scenarios, based on the "Ahuja\_400\_400" instance, which involves 400 targets and 400 weapons.

\begin{table}[!ht]
\centering
\caption{First test case based on Ahuja\_400\_400 - the number of available weapons is equal to the number of targets. Experimental results are presented for the Weapon-Target Assignment (WTA) problem, which involves engaging targets over three consecutive time intervals. The table summarizes the number of weapons available for engagement after each interval, the number of targets remaining to be engaged after each interval, and the number of targets that remained unneutralized at the final stage. For each stage, the average time required to determine a solution and the average value of the objective function are provided.}
\label{tab::DWTA1}
\resizebox{\textwidth}{!}{

\begin{tabular}{llllll}
\toprule
case & average  & standard & average & average & average \\
 &  score & score & time & targets left & weapons left \\
\midrule
SWTA & 2457.72 & 389.21 & 1.43 & 40 & 0 \\
DWTA-A rnd & 13972.00|6278.53|1277.67 & 315.66|305.83|223.58 & 2.34|0.77|0.19 & 279|159|40 & 266|132|0 \\
DWTA-A inc. & 13910.94|6232.38|1268.72 & 293.92|275.83|217.85 & 2.84|0.87|0.20 & 278|159|40 & 265|132|0 \\
DWTA-A dec. & 13996.38|6257.23|1323.07 & 305.60|298.08|222.52 & 2.87|0.87|0.11 & 279|159|41 & 266|132|0 \\
DWTA-A eq. & 13989.83|6320.53|1332.48 & 330.14|291.37|215.31 & 2.81|0.84|0.15 & 279|160|42 & 266|132|0\\
\bottomrule
\end{tabular}
}
\end{table}

It can be observed that, regardless of the scenario regarding changes in probabilities, the application of a "kill-assessment" step between successive stages results in an approximately twofold improvement in the expected value of the objective function. At the same time, the number of targets that remain undestroyed is the same as in the static version of the problem. This is characteristic of situations where the defense does not possess more engagement assets than there are threats. In such cases, the solver prioritizes targets with the highest value, while disregarding those deemed less dangerous.

Table \ref{tab::DWTA2} presents calculations for both static and dynamic Weapon-Target Assignment (WTA) scenarios, based on the "Ahuja\_40\_20" instance. In this specific case, the number of available weapons is twice the number of targets—a favorable situation for defensive systems.
\begin{table}[!ht]
\centering
\caption{Experimental results for the dynamic and static WTA case based on Ahuja\_40\_20 set – number of weapons exceeds number of targets.}
\label{tab::DWTA2}
\begin{tabular}{llllll}
\toprule
 case & average  & standard & average & average & average \\
 &  score & score & time & targets left & weapons left \\
\midrule
SWTA & 19.68 & 39.36 & 0.03 & 0.30 & 0 \\
DWTA-A rnd & 422.20|24.21|0.00 & 79.81|34.89|0.00 & 0.03|0.02|0.01 & 8|0.54|0 & 27|13|0 \\
DWTA-A inc. & 505.48|27.85|0.00 & 94.65|41.96|0.00 & 0.03|0.02|0.01 & 9|0.59|0 & 28|14|0 \\
DWTA-A dec. & 323.33|9.51|0.00 & 85.71|26.88|0.00 & 0.03|0.02|0.00 & 7|0.21|0 & 25|13|0 \\
DWTA-A eq. & 444.29|33.88|0.00 & 100.34|46.83|0.00 & 0.03|0.03|0.01 & 9|0.69|0 & 27|13|0 \\
\bottomrule
\end{tabular}
\end{table}
In a scenario where the number of weapons available is twice the number of targets to be attacked, the expected value—already small—of both the number and the value of surviving targets is reduced to zero across this range of probabilities. In the case of DWTA tasks, virtually all targets are eliminated as early as the second time interval.

Table \ref{tab::DWTA3} was developed based on the Ahuja\_20\_80 instance. In this case, the number of targets exceeds the number of weapons, a situation that is disadvantageous for the defense system.

\begin{table}[!ht]
\centering
\caption{Experimental results for the dynamic and static variants of the Weapon-Target Assignment (WTA) problem, based on the Ahuja\_20\_80 dataset – the number of available weapons is significantly lower than the number of targets.}
\label{tab::DWTA3}
\resizebox{\textwidth}{!}{
\begin{tabular}{llllll}
\toprule
case & average  & standard & average & average & average \\
 &  score & score & time & targets left & weapons left \\
\midrule
SWTA & 3492.41 & 130.41 & 0.02 & 62 & 0 \\
DWTA-A rnd & 4491.43|3962.52|3480.27 & 70.53|124.07|112.86 & 0.03|0.02|0.01 & 73|67|62 & 13|6|0 \\
DWTA-A inc. & 4847.98|4197.69|3478.27 & 60.72|112.73|114.43 & 0.03|0.02|0.01 & 77|70|62 & 17|9|0 \\
DWTA-A dec. & 4335.77|3769.77|3503.88 & 83.90|113.98|113.72 & 0.04|0.02|0.01 & 71|65|62 & 11|4|0 \\
DWTA-A eq. & 4617.35|4109.81|3530.81 & 86.08|98.91|135.12 & 0.03|0.02|0.01 & 74|69|62 & 14|8|0\\
\bottomrule
\end{tabular}
}
\end{table}
When there is an excess of targets, no improvement can be observed—at least within the tested range of probability variations.

Table \ref{tab::DWTA4} was prepared for the Ahuja\_100\_200 instance. In this case, the number of targets exceeds the number of available weapons, representing the most unfavorable scenario for defense systems. Regardless of how target-destruction tasks are assigned across the three time periods, some targets inevitably remain undestroyed. Given the large number of targets that cannot be neutralized due to a shortage of weapons, certain destruction scenarios may yield worse outcomes than those obtained by solving the static version of the problem.

\begin{table}[!ht]
\centering
\caption{Experimental results for the dynamic and static variants of the Weapon-Target Assignment (WTA) problem, based on the Ahuja\_100\_200 dataset – the number of available weapons is significantly lower than the number of targets.}
\label{tab::DWTA4}
\resizebox{\textwidth}{!}{
\begin{tabular}{llllll}
\toprule
case & average  & standard & average & average & average \\
 &  score & score & time & targets left & weapons left \\
\midrule
SWTA & 4883.78 & 250.00 & 0.10 & 110 & 0 \\
DWTA-A rnd & 9299.08|6789.17|4694.09 & 171.46|196.50|190.87 & 0.27|0.12|0.04 & 170|140|110 & 67|33|0 \\
DWTA-A inc. & 9308.42|6797.90|4711.14 & 182.18|211.08|194.35 & 0.31|0.13|0.04 & 170|140|110 & 67|34|0 \\
DWTA-A dec. & 9187.36|6619.06|4747.69 & 149.96|181.24|195.90 & 0.32|0.10|0.03 & 169|138|110 & 66|31|0 \\
DWTA-A eq. & 9241.21|6774.15|4762.28 & 181.25|244.47|222.29 & 0.24|0.10|0.03 & 169|140|111 & 66|33|0 \\
\bottomrule
\end{tabular}
}
\end{table}

In larger instances with a greater number of objectives, it is possible to achieve an improvement in the objective function. For such instances, the solver has more flexibility in selecting the objectives to which it allocates countermeasures.

\begin{table}[!ht]
\centering
\caption{Experimental results for the dynamic and static WTA case based on Ahuja\_200\_100 set – number of weapons exceeds number of targets.}
\label{tab::DWTA5}
\begin{tabular}{llllll}
\toprule
case & average  & standard & average & average & average \\
 &  score & score & time & targets left & weapons left \\
\midrule
SWTA & 62.55 & 60.03 & 0.44 & 1 & 0 \\
DWTA-A rnd & 1674.34|53.67|0.00 & 171.15|43.85|0.00 & 1.19|2.35|0.02 & 39|2|0 & 133|66|0 \\
DWTA-A inc. & 1756.71|61.74|0.00 & 202.02|47.03|0.00 & 0.71|2.95|0.02 & 41|2|0 & 134|67|0 \\
DWTA-A dec. & 1696.16|50.70|0.00 & 179.93|47.16|0.00 & 0.86|0.18|0.02 & 39|1|0 & 133|66|0 \\
DWTA-A eq. & 1719.68|60.52|0.00 & 181.19|49.50|0.00 & 0.45|3.59|0.14 & 40|2|0 & 133|66|0 \\
\bottomrule
\end{tabular}
\end{table}

Table \ref{tab::DWTA5} was developed based on the Ahuja\_200\_100 instance. In this scenario—which is favorable for defense systems—there is a clear numerical superiority of weapons over targets. As was the case with the smaller instance, the kill-assessment process further reduces the already low survivability of the targets. In virtually every instance of the dynamic threat engagement, the threats are eliminated as early as the second operational period of the systems.

\section{Discussion}
\label{sec:Discussion}
The classical approach to solving the WTA problem involves linearizing the problem and solving it as a Mixed-Integer Programming (MIP) task. The paper by Andersen \cite{andersen} presents three approaches for linearizing the WTA problem and transforming it into an MIP formulation. Efficient solvers such as CPLEX or Gurobi can be applied to problems of this class.
The best WTA linearization approach presented in \cite{andersen} ensures an exact linearization of the WTA objective function while simultaneously strengthening the linear formulation through the use of a set of convex, piecewise-linear under-approximating functions.
The paper introduces an exact algorithm that avoids the issue of an excessive increase in the size of the corresponding linear programming models.
This algorithm corrects the value of each solution obtained by the solver, adjusting it to the true value derived from the original objective function.
Experiments conducted in \cite{andersen} confirmed the improved efficiency of the presented method compared to the classical nonlinear problem formulation.
The linearization-based approach, utilizing so-called probability chains, enabled the effective determination of optimal solutions for problems involving up to 350 weapon types, under the additional assumption that at most one weapon unit of a given type could be assigned to a single target. For larger-scale problems, finding an optimal solution becomes a significant challenge; in such cases, one could rely on the heuristic algorithms presented in Section \ref{sec:RelatedWorks}.

To determine an optimal or near-optimal solution, it was necessary to rely on Mixed-Integer Nonlinear Programming (MINLP) models. To find the solution efficiently, an approach based on objective function linearization and column generation was employed. This method reduces the solution time by several orders of magnitude compared to classical methods. The paper presents three possible approaches to solving the Weapon-Target Assignment (WTA) problem. The first involves linearizing the objective function, which is then used in the column generation (CG) subproblem during the pricing procedure. The second method does not require objective function linearization \cite{chen}, as the subproblem is solved by enumerating possible weapon-to-target assignment configurations. The third approach, based on the work of Bertsimas \cite{bertsimas}, features a simplified pricing step that does not require solving an Integer Programming (IP) problem.

The analyses presented in the article demonstrate the high performance of the method proposed by Bertsimas et al. Nevertheless, the method based on the linearized problem proves to be highly efficient in finding solutions. An advantage of this method over the Bertsimas approach is the ease with which the model can be extended by adding constraints to the WTA problem; with the concept presented in \cite{bertsimas}, this process is far more complex.
It is also worth noting that the method based on the linearized problem scales linearly and does not entail significant memory complexity (the subproblems are small IP tasks solved by the solver in, at most, a few seconds).

\subsection{Possible approaches to DWTA problem modeling}

Receding horizon optimization is a control method that addresses the Weapon-Target Assignment (WTA) optimization problem over a finite time horizon—specifically, within defined time intervals during which threats can be engaged. This approach allows for the dynamic adjustment of the allocation of defensive assets against incoming threats across successive time intervals, starting from the present moment. Initially, a solution is determined for the entire set of intervals; the control window then shifts forward, and the process is repeated. This technique is also widely known as Model Predictive Control (MPC).

MPC can be applied to solve Dynamic Weapon-Target Assignment (DWTA) problems. In virtually all the cases presented, employing this method enhances the effectiveness of threat engagement compared to static approaches. This work explores several variants regarding how the probability of a target being destroyed by a specific weapon might change. While many other variants could exist, discussing them in detail would require referencing specific types of defensive assets and aerial threats, which falls outside the scope of this work.
The results presented demonstrate the utility of the receding horizon method when applied to aerial threat scenarios using a WTA algorithm.

In this paper, we consider the deterministic, dynamic Weapon-Target Assignment (WTA) problem and propose an approach based on the Receding Horizon method, which is fully justified in this deterministic context. Primarily, this technique can be applied if the optimal control problem for each time horizon can be solved within a short timeframe—a condition that is well-met in the case of our control problem. Furthermore, for a deterministic dynamic task, the solution obtained using the RH technique is identical to the solution derived via the Bellman method.

The Receding Horizon approach offers a certain degree of robustness against system uncertainties,
its deterministic formulation generally makes it inherently unsuitable for systematically dealing with
uncertainties. We could treat the optimal control problem—which we solve as part of the receding horizon (RH) strategy—as a problem involving an uncertain system; however, we then encounter a difficulty, namely that the technique commonly used to solve such a problem is one based on Bellman equations.As a result, extending the presented approach to the stochastic case does not appear trivial. Research aimed at developing a stochastic version of the RH approach to the dynamic WTA problem could draw upon findings discussed in the review papers \cite{Mesbah2016} and \cite{Farina2016}.

\section{Conclusions}
\label{sec:Conclusions}
This article presents algorithms for solving Weapon-Target Assignment (WTA) problems, based on the column generation technique. Methods for determining optimal solutions for both static and dynamic versions of the problem are described. The version of the Dynamic Weapon-Target Assignment (DWTA) problem presented here relies on the receding horizon technique. Existing literature addresses the dynamic WTA problem primarily through Bellman’s dynamic programming equations; however, this approach is only feasible for small-scale problems. In contrast, the method proposed in this work—utilizing the receding horizon technique—is scalable, enabling the solution of large-scale dynamic WTA problems.

This paper considers a deterministic variant of the dynamic WTA problem; developing an effective method for solving the stochastic dynamic WTA problem remains an open issue.

\bibliographystyle{unsrt}  
\bibliography{references}  


\end{document}